\documentclass[mnsc,nonblindrev]{informs4}

\OneAndAHalfSpacedXII % Current default line spacing
\usepackage{natbib}
 \bibpunct[, ]{(}{)}{,}{a}{}{,}%
 \def\bibfont{\small}%

\def\N{\mathbb{N}}
\def\R{\mathbb{R}}

\def \E {\mathbb{E}}
\def \P {\mathbb{P}}

\def \UC {\mathcal{U}}
\def \-> {\rightarrow}
\def \ind {\mathbf{1}}

\def \xb {\mathbf{x}}
\def \yb {\mathbf{y}}
\def \Xb {\mathbf{X}}

\def \eqd {\,{\buildrel d \over =}\,}

\def \argmin {\operatorname{argmin}}

\def \Var {\mbox{\bf Var}}
\def \OT {\Tilde{O}}
\def \XC {\mathcal{X}}
\def \XbT {\Tilde{\mathbf{X}}}
\def \xbT {\Tilde{\mathbf{x}}}
\def \ybT {\Tilde{\mathbf{y}}}
\def \nubar {\overline{\nu}}

\def \Db {\mathbf{D}}
\def\DeltaQ{\Delta Q^r}

\usepackage{enumerate}
\usepackage{enumitem}
\usepackage{tabularx,booktabs,mathrsfs,setspace}%
\usepackage{mathtools}
\usepackage{multirow}
\usepackage{rotating}

\TheoremsNumberedThrough     % Preferred (Theorem 1, Lemma 1, Theorem 2)
\EquationsNumberedThrough    % Default: (1), (2), ...
\begin{document}
%%%%%%%%%%%%%%%%

%\setcounter{page}{1}
%\DOI{10.1287/xxxx.0000.0000}%

\RUNAUTHOR{Moradi et al.}

\RUNTITLE{Asymptotic Optimality of Inventory Projection in Dual-Sourcing Systems}

\TITLE{Asymptotic Optimality of Inventory Projection in Dual-Sourcing Systems}

\ARTICLEAUTHORS{
\AUTHOR{Poulad Moradi}
\AFF{School of Industrial Engineering, Eindhoven University of Technology, Eindhoven, the Netherlands, \EMAIL{p.moradi.shahmansouri@tue.nl}}
\AUTHOR{Melvin Drent}
\AFF{Department of Information Systems and Operations Management, Tilburg University, Tilburg, the Netherlands, \EMAIL{m.drent@tilburguniversity.edu}}
\AUTHOR{Joachim Arts}
\AFF{Luxembourg Centre for Logistics and Supply Chain Management, University of Luxembourg, Luxembourg-City, Luxembourg,  \EMAIL{joachim.arts@uni.lu}}
\AUTHOR{Collin Drent}
\AFF{School of Industrial Engineering, Eindhoven University of Technology, Eindhoven, the Netherlands, \EMAIL{c.drent@tue.nl}}
} 

\ABSTRACT{%
We consider a single-echelon inventory system under periodic review with two suppliers facing stochastic demand, where excess demand is backlogged. The expedited supplier has a shorter lead time than the regular supplier but charges a higher unit price. We introduce the Projected Expedited Inventory Position (PEIP) policy. We show that this policy is asymptotically optimal in two regimes of practical interest: (i) when the lead time of the regular supplier becomes large and (ii) when both the shortage cost and the cost premium for expedited units become large, with their ratio held constant. We show through an extensive numerical investigation that the PEIP policy outperforms the current best performing heuristic policies in literature.
}%

%\KEYWORDS{deterministic inventory theory; infinite linear programming duality; 
%  existence of optimal policies; semi-Markov decision process; cyclic schedule}

\KEYWORDS{dual-sourcing; inventory; projection policies; asymptotic optimality; long lead times; high service}
%\HISTORY{submitted january 2020}

\maketitle
%%%%%%%%%%%%%%%%%%%%%%%%%%%%%%%%%%%%%%%%%%%%%%%%%%%%%%%%%%%%%%%%%%%%%%

% Samples of sectioning (and labeling) in MSOM
% NOTE: (1) \section and \subsection do NOT end with a period
%       (2) \subsubsection and lower need end punctuation
%       (3) capitalization is as shown (title style).
%
%\section{Introduction.}\label{intro} %%1.
%\subsection{Duality and the Classical EOQ Problem.}\label{class-EOQ} %% 1.1.
%\subsection{Outline.}\label{outline1} %% 1.2.
%\subsubsection{Cyclic Schedules for the General Deterministic SMDP.}
%  \label{cyclic-schedules} %% 1.2.1
%\section{Problem Description.}\label{problemdescription} %% 2.

\section{Introduction}
\label{sec:intro}
Dual-sourcing is a common and effective supply chain management practice employed by companies that operate both globally and
domestically. It generally refers to situations where companies can replenish their inventories from a regular supplier as well as from an expedited supplier; 
the latter offering a shorter lead time than the former but at the expense of a cost premium.
Even though the expedited supplier charges an additional cost, companies can leverage its responsiveness when
inventories deplete precipitously.

Examples of companies that have utilized a dual-sourcing strategy are manifold. \cite{allon2010global} describe a
\$10 billion high-tech U.S. company with suppliers in China and Mexico. 
The Chinese supplier charges lower costs than Mexico but shipping its products to the U.S. takes considerably longer.
\cite{rao2000caterpillar} describe a similar situation at Caterpillar; they too use a more expensive expedited supplier whenever a shorter lead time is desired. 
While both examples concern two distinct suppliers, dual-sourcing may also refer to situations where companies can replenish their inventories from one supplier using two distinct transport modes. 
For example, \cite{threatte2002tactical} describe a situation in which Polaroid Corporation relies on ocean transport to ship products from a supplier in Asia to the U.S., but occasionally uses expedited air transport when inventory levels are critically low.
The choice between transportation modes may also be used to jointly reduce costs and the overall carbon footprint of a supply chain \citep{Drent2025}.

In this paper, we treat a setting similar to the ones described above. In particular, we consider a single-item, single-echelon dual-sourcing inventory system under periodic review facing stochastic demand, where excess demand is backlogged. 
We assume that demands across periods are independently and identically distributed (i.i.d.).
The expedited supplier has a shorter lead time than the regular supplier but charges a higher unit price. 
There are no fixed costs associated with ordering at either supplier. 
At the end of each period, linear holding and backorder costs are charged for leftover inventory and backlogged demand, respectively. We seek to minimize the cost-rate, which consists of ordering, holding, and backorder costs.  

The dual-sourcing system introduced above has been studied for almost 65 years now, starting with the contributions of \cite{barankin1961delivery} and \cite{neuts1964}. 
They show that the optimal policy exhibits a simple order-up-to structure when the lead times of the regular and expedited supplier are one and zero periods, respectively. \cite{fukuda1964optimal} later extends this result to the case where the difference between the lead times of both suppliers is one period, irrespective of their individual lead times. 
For lead time differences larger than one period, \cite{whittemore1977optimal} and later \cite{feng2006base} show that the structure of the optimal policy is complex, and depends in general on the entire vector of outstanding orders. 
Although the optimal policy for general lead time differences can be computed in principle by stochastic dynamic programming, this is not practical as the state space of the dynamic program grows exponentially in the lead time difference; that is, it suffers from the curse of dimensionality. 
Most researchers have therefore focused on devising relatively simple but non-optimal replenishment policies (we will provide a thorough review later in Section \ref{sec:literatureReview}). 

Even though some heuristic policies have been shown to perform numerically close to optimal in at least certain regimes, theoretical justifications for such useful asymptotic optimality performance guarantees are mostly lacking; the only exception being the Tailored Base-Surge (TBS) policy. 
This policy orders a fixed amount at the regular supplier every period and uses an order-up-to rule for placing expedited orders. 
The TBS policy was first proposed by \cite{allon2010global}, though \cite{rosenshine1976analysis} and \cite{janssen1999two} have studied a closely related policy before.
Numerical analysis by \cite{klosterhalfen2011comparisonDIP} and \cite{janakiraman2015analysis} indicates that the optimality gap of the TBS policy closes as the lead time difference grows large. \cite{janakiraman2015analysis} also establish a near-optimality result under a demand structure consisting of a base component and an infrequent surge component

\cite{xin2018asymptotic} provide theoretical justification for these numerical observations by proving that a simple TBS policy is asymptotically optimal as the lead time of the regular supplier grows large, with the lead time of the expedited supplier held fixed for a general demand distribution.
More recently, \cite{Chen2026} show that the optimality gap of TBS decays exponentially with the regular supplier's lead time.
The intuition behind this result is that as the lead time difference grows, the amount of randomness between when a regular order is placed and when this order arrives grows so large that a constant order, which ignores the entire vector of outstanding orders, performs nearly optimal.
To our knowledge, this is the only heuristic policy for dual-sourcing inventory systems under periodic review with general lead time differences that has an asymptotic  performance guarantee. 
Unfortunately, the TBS policy performs poorly outside this asymptotic regime, particularly when either the backorder cost becomes large or when
the unit price charged by the expedited supplier becomes large. When both cost parameters become large simultaneously,
the performance of the TBS drops even precipitously \citep{klosterhalfen2011comparisonDIP}. Yet, it is exactly this setting that
is prevalent in practice, and it therefore remains an open research problem to determine control policies that are provably asymptotically optimal in this regime.

In this paper, we fill this gap and propose such a new policy, which we call the Projected Expedited Inventory Position (PEIP) policy. 
The PEIP policy keeps track of the expedited inventory position (net inventory level plus outstanding orders arriving within the expedited lead time) and places expedited orders such that this inventory position reaches a certain target level. 
Regular orders entering the information horizon for placing expedited orders may cause the expedited inventory position to exceed its corresponding target level. 
This implies that the expedited inventory position in expectation exceeds the target level. 
Under the PEIP policy, regular orders are therefore placed such that the expected expedited inventory position at the time of the regular order entering the expedited inventory position is kept at a target level. 
In doing so, the PEIP policy leverages all the information contained in the pipeline -- as does the optimal policy.
This is in contrast with existing heuristic policies that either aggregate all outstanding orders or ignore them entirely and place a constant order every period regardless. 
We further provide a separability result that enables efficient optimization of the PEIP policy. 
Thus, despite exploiting substantially richer pipeline information than existing heuristics, the PEIP policy remains computationally tractable and retains an intuitive appeal for practitioners. Moreover, it is reassuring to note that, in the special case of consecutive lead times, the PEIP policy reduces to a policy characterized by two order-up-to levels, which is known to be optimal in that setting \citep{fukuda1964optimal}. 
We then study the performance of the PEIP policy in two asymptotic regimes of interest. In the first regime, the lead times of the regular supplier grow large. We assume that demand per period has finite variance and show that in this regime, the PEIP policy matches the performance of the TBS policy. This implies in particular that the PEIP policy is also asymptotically optimal as the regular lead time grows large. Our approach begins by establishing that the state process admits an stationary probability measure under the PEIP policy. We then study a coupling with the TBS policy and prove a contraction inequality that shows that the regular orders generated by PEIP converge to a constant as the lead time grows. 

In the second asymptotic regime, both the shortage cost and the price for expedited units become large, with their ratio held constant. We call this the high service regime and show that, under some mild conditions on the demand distribution, the relative difference between the cost-rate of the best PEIP policy and the optimal policy converges to zero. Our approach relies on establishing lower and upper bounds on the cost-rates of the optimal policy and the best PEIP policy, respectively, in terms of the cost-rates of canonical inventory systems with only one supplier. We then show that the ratio of these bounds converge to 1 in this high service regime where both the shortage cost and the price for expedited units grow large.

Our research is useful for practice as it directly leads to an easily implementable heuristic that has asymptotic performance guarantees in systems with long lead times or high backlog penalties and an expensive expedited source.  Both features are prevalent in practice, sometimes in conjunction. The capital goods industry provides one prominent example, as shortages of critical parts can be very costly and the expedited source may be an expensive 3D printer \citep[see, e.g.,][]{westerweel2018printing}. The lead-time difference between local 3D printing and traditional manufacturing followed by shipment across the globe is also substantial.

The main contributions of this paper are fourfold:
\begin{enumerate}
\item We propose the PEIP policy for periodic review dual-sourcing inventory systems with general lead time differences, and we present an efficient
optimization scheme to find its optimal policy parameters. This optimization scheme hinges on a separability result that allows us to decompose the optimization into two simple one-dimensional optimization problems.
\item We prove that the PEIP policy is asymptotically optimal when the lead time of the regular supplier becomes large. Our proof technique is new and relies on first establishing the existence of a unique stationary probability measure under PEIP and then exploiting a contraction inequality that arises out of a coupling with the TBS policy. We also explain how this proof technique extends to PIL policies for the canonical lost sales inventory model, a special case of the dual-sourcing model considered here \citep{sheopuri2010newPolicies,vanJaarsveld2024}, and may provide a template for other projection-based inventory policies.
\item We prove that the PEIP policy is asymptotically optimal when both the shortage cost and the cost premium  for  expedited  units  become large,  with their ratio held constant and call this the high service regime.
Our proof technique can serve as a template that is readily applicable to other heuristic policies. We illustrate this by showing that several existing heuristic policies are also asymptotically optimal in this regime.
\item Through an extensive numerical investigation, we conclude that the PEIP policy outperforms the best performing heuristics in literature.
\end{enumerate}

The remainder of this paper is organized as follows. In the next section, we position our research relative to the existing literature. After that, in Section \ref{sec:description}, we formally define the decision problem and introduce the notation that we use throughout this paper. 
Section \ref{sec:polPolicy} presents the PEIP policy, including analytical results that allow for an efficient optimization procedure to find the optimal policy parameters. 
We proceed with establishing our asymptotic optimality result for long regular lead times in Section~\ref{sec:longleadtime} and for the high service regime in Section \ref{sec:asymptoticOptimality}. 
After that, we report on a large numerical study in Section \ref{sec:numericalStudy}, and provide concluding remarks in Section \ref{sec:concludingRemarks}.

\section{Positioning in the Literature}
\label{sec:literatureReview}
Dual-sourcing systems operating in discrete time as well as their variants (e.g., continuous review models, multi-echelon systems, etc.) have been widely studied. 
For a broad overview of this rich field, we refer the interested reader to the surveys of \cite{minner2003multiple} and  \cite{svoboda2019review}, as well as to the recent studies of, e.g.,  
\cite{gong2014dynamic}, 
\cite{boute2015global}, \cite{arts2016repairable}, \cite{sapra2017dual},
\cite{song2017optimal}, \cite{drent2020expediting}, \cite{deng2026dual} and the references therein. Here we confine ourselves to contributions most relevant to our research, a large part of which revolves around heuristic policies.

\citeauthor{whittemore1977optimal} show already in 1977 that optimal policies for the periodic review dual-sourcing inventory systems we typically face in practice (i.e. with lead time differences larger than one) are complex and difficult to compute. 
Yet it lasted until the beginning of this century before researchers began to devise various effective heuristic policies.
\cite{scheller2007effective} propose the Single-Index (SI) policy that uses order-up-to rules for order placement at the regular and expedited supplier based on only the regular inventory position (inventory level plus all outstanding orders).
By additionally tracking the expedited inventory position that includes only orders that will arrive within the expedited lead time, \cite{veeraraghavan2008now} introduce the Dual-Index (DI) policy.
Regular and expedited orders are then placed according to order-up-to rules based on their respective inventory positions. 
They propose a simulation-based optimization and show that this is easily implementable and often results in near-optimal solutions in numerical examples. 

Since then, the DI policy has been extended in various useful ways. \cite{sheopuri2010newPolicies} analyze the more general class of policies that use an order-up-to rule for placing expedited orders.  
\cite{arts2011efficient} incorporate stochastic lead times and provide an approximate evaluation method so that optimization can be done without simulation. 
\cite{sun2018capped} introduce the Capped Dual-Index (CDI) 
policy in which the DI policy is accompanied with a cap that smooths regular orders. 
%That is, the CDI policy operates as a DI policy, except that the regular order size cannot
%exceed its cap. 
Even though the CDI policy has been developed in a robust optimization setting under the worst-case cost performance criterion, numerical analysis indicates that it also performs well in terms of minimizing the cost-rate in a stochastic setting.  

\cite{sheopuri2010newPolicies} show that the inventory system with lost sales is a special case of the dual-sourcing inventory system treated here given that the expedited supplier follows an order-up-to rule. 
The intuition behind this is that the orders placed at the expedited supplier can be thought of as demand lost from the regular supplier, with the cost premium of the expedited supplier being the lost sales penalty cost. 
Inspired by this relation, \cite{sheopuri2010newPolicies} propose the Vector Base-Stock (VBS) policy that uses an order-up-to rule for the expedited supplier, and a vector base-stock policy for the regular supplier.
\cite{hua2015structural} later generalize the VBS policy to the Best Weighted Bounds (BWB) policy. In addition, by using properties of $L^{\natural}$-convexity, they show that the optimal expedited orders are more sensitive to the older (i.e. soon-to-arrive) outstanding orders whereas the optimal regular orders are more sensitive to the younger outstanding orders.
Extensive numerical investigations indicate that the BWB policy as well as the CDI policy described earlier are the best performing heuristic policies in the current literature \citep{sun2018capped}, and can even outperform advanced deep reinforcement learning algorithms \citep{gijsbrechts2019drl}.    

%hier beneden nog aanpassen; ik bedoel dat alle policies een overshoot hebben, maar dat wij intuitief daarop sturen.
While the PEIP policy is similar to the aforementioned heuristic policies in that it also assumes an order-up-to rule for the expedited supplier, it differs significantly in how it operates the regular supplier. 
First, the PEIP policy projects the expedited inventory position in deciding upon regular orders. 
As such, it controls the expected expedited inventory position directly taking into account explicitly the excess by which the order-up-to level of the expedited supplier overshoots; existing heuristics do so only indirectly. 
Second, similar to the optimal policy, the PEIP policy leverages all available state information in deciding upon order placement at the regular supplier. 
Existing heuristics either aggregate all available information, or ignore it entirely and place a fixed amount every period regardless.

Our research also contributes to the growing body of literature that establishes asymptotic performance guarantees of heuristic policies for inventory systems whose optimal policies are difficult to analyze if not intractable. This stream was reviewed by \cite{goldberg2019survey}. Within this stream, many papers study the single-echelon inventory systems with lost sales.  
For such systems, \cite{huh2009asymptotic} and \cite{bijvank2014robust} show that certain order-up-to policies are asymptotically optimal as the lost sales penalty cost approaches infinity while \cite{goldberg2016asymptotic} and \cite{xin2016asymptotic} show that a simple constant order policy is asymptotically optimal as the lead time approaches infinity.
Both heuristic policies perform however poorly outside their respective asymptotic regimes. 
This has recently led to the development of two heuristic policies that have asymptotic performance guarantees and additionally perform well in both regimes: the Capped Base-Stock (CBS) policy proposed by \cite{xin2020approx,xin2020capped}, and the Projected Inventory Level (PIL) policy proposed by \cite{vanJaarsveld2024}. 

The CBS policy can be interpreted as a hybrid between a constant order policy and an order-up-to policy, and it thus enjoys the superior performance of both policies in their respective regimes. 
The PIL policy places orders  such  that  the  expected  inventory level at the time of arrival of an order is raised to a target level. 
\cite{vanJaarsveld2024} show that the PIL policy is asymptotically optimal for large lost sales penalty costs as well as for large lead times, though the latter only under the assumption that demand is exponentially distributed. The need for this assumption has recently been relaxed by \cite{Moradi2025}.
The PEIP policy is related to the PIL policy:
While a PIL policy places orders by projecting the expected inventory level, a PEIP policy places regular orders by projecting the expected expedited inventory position (while retaining an order-up-to rule for the expedited supplier). Projection type policies also have desirable asymptotic optimality properties in perishable inventory systems \citep{bu2026managing} and multi-echelon inventory systems with lost-sales \citep{Gong2026}.

Unlike inventory systems with lost sales, there currently exists only one heuristic policy for periodic review dual-sourcing inventory systems that has asymptotic performance guarantees.
Indeed, a simple TBS policy is asymptotically optimal as the lead time difference between the regular supplier and expedited supplier grows large \citep{xin2018asymptotic} and the optimality gap decays exponentially with the lead time difference \citep{Chen2026}. The fact that the regular order does not depend on the current system state under a TBS policy is a generic feature of open-loop control. \citet{Bai2023} show that such open-loop controls are asymptotically optimal in a broader class of MDPs, including dual-sourcing systems.  \cite{bai2026asymptotic} extend this result to long random lead times. We prove the asymptotic optimality of PEIP policies in the long lead time regime by first establishing that a unique invariant probability measure exists under PEIP and then exploiting a contraction inequality in a coupled system with the TBS. This contraction shows that the PEIP converges to a TBS when regular lead times become large.

In addition to the long lead time regime, this paper considers a new asymptotic regime in which the unit cost of the expedited supplier and the shortage cost approach infinity simultaneously. We call this regime the high service regime 
and show that the best PEIP policy is asymptotically optimal in this regime. 
Our main proof exploits the earlier described relation between lost sales inventory systems and dual-sourcing inventory systems; this allows us to draw on the asymptotic results for lost sales inventory systems of \cite{huh2009asymptotic}. 
Our proof technique can serve as a template that is readily applicable to existing as well as future heuristic policies. We illustrate this by showing that the CDI, DI, and SI policy are also asymptotically optimal in the high service regime. 

%In addition, we show numerically that a PEIP policy converges to a TBS policy when the lead time difference grows large. As such, 
%As such, our work essentially complements the work of \cite{xin2018asymptotic}.
%Indeed, as \citet[p. 445]{xin2018asymptotic} remark,``it would be interesting to identify other more sophisticated algorithms that perform better for small-to-moderate lead times''.
We end this review by remarking that recent advances on asymptotic optimality results for inventory systems that have complex optimal policies stretch beyond the lost sales and dual-sourcing inventory systems described above. 
See, for instance, \cite{Rong2017heuristics} for such results in the context of distribution systems consisting of one central warehouse and multiple local warehouses, \cite{janakiraman2018simple} on inventory systems with multiple products that share production capacities, and \cite{fu2019managing} on inventory systems for perishable products.  
We refer the interested reader to \cite{goldberg2019survey} for a recent and comprehensive review of this growing field.

\section{Notation and Problem Formulation}
\label{sec:description}
This section provides a formal problem description that mostly follows the notational conventions of \cite{sheopuri2010newPolicies}.
We consider a single-item, single-echelon inventory system under periodic review facing stochastic demand, where
excess demand is backlogged. 
Products can be ordered each period both at a regular supplier and at an expedited supplier.
Let $c^r$ ($l^r$) and $c^e$ ($l^e$) denote the unit prices (lead times) from the regular and expedited supplier,
respectively. 
We assume that $c^e>c^r$ and $l^e<l^r$, that is, the expedited supplier offers a higher unit price, but its lead
time is shorter. 
Notice that the present inventory problem reduces to a standard inventory problem with only one supplier if these
assumptions do not hold.
Both lead times are assumed to be deterministic as well as a non-negative integer multiple of the review period. 
For notational convenience, we also define the lead time difference $l=l^r-l^e \geq 1$. 

Let $D_t$ denote the random demand in period $t$, with $t$ indexed forward in time starting from 0, i.e.
$t\in\{0,1,\ldots\}$.
Demand across periods is a sequence of non-negative i.i.d. random variables with finite first moment $\mu:=\E[D_0] \in (0,\infty)$.
Let $\sigma^2:=\Var[D_0]$. We will make mild additional assumptions on the demand distribution in Section~\ref{sec:longleadtime} and Section~\ref{sec:asymptoticOptimality}.
Whenever we drop the period index of $D_t$, we refer to the generic random variable with the same distribution as $D_t$. The cumulative demand over $n$ periods is denoted as $\Db^n$.
Let $I_t$ denote the net inventory level (on-hand inventory minus backlog) at the beginning of period $t$ before orders and demand arrive. 
Let $Q^e_t$ $(Q^r_t)$ denote the order size placed with the expedited (regular) supplier in period $t$.
The orders that arrive in period $t$ are thus $Q^e_{t-l^e}$ and $Q^r_{t-l^r}$, and we can consequently write the following recursion for the inventory level:
\[
I_{t+1} = I_t - D_t + Q^e_{t-l^e} + Q^r_{t-l^r}.
\]
Let $J_t$ denote the net inventory level after receiving orders, i.e., $J_t:= I_t+Q^e_{t-l^e}+ Q^r_{t-l^r}$.
At the end of each period $t$, the on-hand inventory carried over to the next period is given by $(J_t-D_t)^+$ and the backlog by $(D_t-J_t)^+$.
Here we use the convention $x^+=\max(x,0)$.
Costs are levied as follows.
The per unit holding cost for each unit of on-hand inventory at the end of period $t \in \N_0$ is $h>0$, and the per unit backlog penalty cost is $p>0$.
% Per unit in on-hand inventory $(I_t+Q^e_{t-l^e}+ Q^r_{t-l^r}-D_t)^+$ carried over to the next period, a cost $h>0$ is charged. 
% Similarly, a cost $p>0$ is charged for each unit in backlog $(D_t-I_t-Q^e_{t-l^e}- Q^r_{t-l^r})^+$ carried over to the next period.

Observe that if the unit price charged by the expedited supplier $c^e$ is higher than the total backlog cost
incurred on a unit over $l$ periods, i.e. $l\cdot p$, then there is no cost benefit to be derived from ordering any
units from the expedited supplier over the regular supplier. In fact, as shown in \cite{sheopuri2010newPolicies},
the present inventory problem reduces to a standard single supplier inventory problem if $c^e \geq p \cdot l$.
We therefore assume throughout the remainder that $c^e< p \cdot l$. We revisit this assumption in Section \ref{sec:asymptoticOptimality} when we will study a high service asymptotic regime.

Let the state of the system at the beginning of period $t$ before already placed replenishment orders arrive be given by $\Xb_t=(I_t,Q^r_{t-l^r},\ldots,Q^r_{t-1}, Q^e_{t-l^e},\ldots,Q^e_{t-1})$.
The sequence of events in any period $t$ is as follows:
\begin{enumerate}
\item The size of the order placed with the expedited supplier $Q^e_t$ is determined based on the current state of the system $\Xb_t$.
\item The size of the order placed with the regular supplier $Q^r_t$ is determined based on the current state of the system $\Xb_t$ and the order just placed with the expedited supplier $Q^e_t$.
\item Products that were ordered $l^e(l^r)$ periods ago at the expedited (regular) supplier are received and added to the on-hand inventory, i.e. order $Q^e_{t-l^e}(Q^r_{t-l^r})$.
\item The demand $D_t$ is realized, and holding or shortage costs are incurred, in case that $(J_t-D_t)^+>0$ or $(D_t-J_t)^+>0$, respectively.
\end{enumerate}  
The total cost incurred in period $t$ is denoted by $C_t$, i.e. $C_t = c^e\cdot Q^e_t + c^r \cdot Q^r_t + h \cdot (J_t-D_t)^+ + p\cdot (D_t-J_t)^+$. 
An admissible policy $u$ is a rule that determines the size of the orders placed with the regular and expedited supplier based only on the historical information available in that period as outlined in the sequence of events above.
All such admissible policies $u$ are contained in the set $\UC$.
Since the total cost incurred in a period $t$ depends on the policy $u$, we henceforth write $C^u_t$ to denote this dependence explicitly. 
Our performance measure of interest is the cost-rate
\[
    C^u = \limsup_{T\rightarrow \infty} \frac{1}{T} \sum\nolimits_{t=1}^T  \E \left[C^u_t \right].
\]

We assume throughout the remainder that $\Xb_0=\mathbf{0}$ unless explicitly stated otherwise.
That is, the initial inventory level equals zero, and there are initially no outstanding expedited or regular orders due to arrive within their corresponding lead times.
This assumption simplifies our analysis later;
it is however not strictly necessary as we are interested in the cost-rate, which is not affected by any specific choice for the initial state.
%Assumption now here as the recursions we previously stated require this assumption too.

We end this section with two remarks related to the optimal policy $u^*=\argmin_{u\in\UC}C^u$. First, as shown by \cite{sheopuri2010newPolicies}, $u^*$ depends on $c_e$ and $c_r$ only through the price premium of the expedited supplier, i.e. $c_e-c_r$. This result follows intuitively from the observation that any policy that has a finite cost-rate will incur at least a cost of $c^r\cdot \mu$.
We will therefore assume without loss of generality that $c_r=0$. 
Second, even for small lead time differences, the structure of $u^*$ remains poorly understood and depends in
general on the entire vector of outstanding orders. 
The policy we propose in the next section also takes into account this entire vector when deciding on the size of
the orders placed with both suppliers. 
However, different from $u^*$, our policy does allow for tractable analysis and efficient optimization, even
for large-sized problems often encountered in practice.    

\section{The Projected Expedited Inventory Position Policy}
\label{sec:polPolicy}
In this section, we present the Projected Expedited Inventory Position (PEIP) policy for operating dual-sourcing inventory systems.  
We first formally define the rules under which both suppliers operate, and we subsequently present analytical results that allow for an efficient solution procedure to find the optimal policy parameters. 

\subsection{Order Placement Rules}
Under a PEIP policy, the ordering rule for the expedited supplier is assumed to be an order-up-to rule that operates as follows.
As a notational convenience throughout, we let $A[t_1,t_2]=\sum_{i=t_1}^{t_2} A_i$ for any time indexed variable $A_t$ and $t_1,t_2\in\N_0$ with $t_1\leq t_2$.
At the beginning of each period $t$ before orders $Q^r_{t-l^r}$ and $Q^e_{t-l^e}$ arrive, we review the expedited inventory position
\begin{equation*}\label{eq:ipe}
IP_t^e=I_t + Q^e[t-l^e,t-1] + Q^r[t-l^r,t-l],
\end{equation*}
based on the reduced state $\Xb^e_t=(I_t,Q^r_{t-l^r},\ldots,Q^r_{t-l}, Q^e_{t-l^e},\ldots,Q^e_{t-1})$.
If necessary, we place an expedited order $Q^e_t$ to raise $IP_t^e$ to its target level $S^e$, i.e. $Q^e_t=(S^e-IP_t^e)^+$.
Note that regular orders entering the information horizon for placing expedited orders may cause the expedited inventory position to exceed this target level. 
That is, immediately after expedited orders have been placed, there is an overshoot defined as $O_t = IP^e_t + Q^e_t - S^e = (IP^e_t - S^e)^+$. 
This implies that the expedited inventory position following expedited ordering is in expectation greater than or equal to the target level. 
In fact, strictly greater unless products are sourced exclusively from the expedited supplier, which we have argued is sub-optimal under the assumptions we make.
The system dynamics can be expressed in terms of the regular order process, $\{Q^r_t\}_{t \in \N_0}$, for a given expediting order-up-to level $S^e$:
\begin{alignat}{2}
    IP_{t+1}^e &=S^e+O_t+Q^r_{t-l+1}-D_t, \label{eq:Dynamics_IP}\\
    O_{t+1} &=(O_t+Q^r_{t-l+1}-D_t)^+, \label{eq:Dynamics_O}\\
    Q^e_{t+1} &=(D_t-O_t-Q^r_{t-l+1})^+, \label{eq:Dynamics_Qe}\\
    J_{t+l^e} &= S^e + O_t - D[t,t+l^e-1],\label{eq:Dynamics_J}\\
    I_{t+l^e+1} &= S^e + O_t - D[t,t+l^e]. \label{eq:Dynamics_I}
\end{alignat}

The PEIP policy aims to keep the expedited inventory position after ordering (and consequently the overshoot) at a fixed level. Although the expedited inventory position cannot actually
be kept at a fixed level, it can be kept at a fixed level in expectation.
More specifically, the PEIP places regular orders such that the expected expedited inventory position at the time those orders enter the information horizon for expedited order placement is raised to a target level.
We henceforth refer to this target level as the projected expedited inventory position level, denoted $\xi$. 

Before we provide a formal description of this ordering rule, we first state the following result related to the minimal required information for optimal regular order placement given that the expedited supplier follows an order-up-to rule.
Let $\Xb^r_t:=(O_t, Q^r_{t-l+1}, \dots,Q^r_{t-1})$.
The information contained in $\Xb^r_t$ is sufficient to determine the expected expedited inventory position in the period that the regular order placed now enters the expedited inventory position. 
\begin{lemma}{\citep[Lemma 4.1]{sheopuri2010newPolicies}} 
Under the assumption that the expedited supplier follows an order-up-to rule, the optimal regular order in period $t$, $Q^r_t$, depends on $\Xb_t$ and $Q^e_t$ only through the $l$-dimensional state description
$\Xb^r_t$. 
\end{lemma}

Observe from Equation~\eqref{eq:Dynamics_IP} that the inventory position in period $t+l$, i.e. $IP^e_{t+l}$, includes the regular order $Q^r_t$ that we will place in the current period $t$.
Hence, $Q^e_{t+l}$ also depends on $Q^r_t$.
We write $\E[IP^e_{t+l} (Q^r_t) +Q^e_{t+l} (Q^r_t) \vert \Xb^r_t]$ to denote the conditional expected expedited inventory position $l$ periods into the future in period $t$ after placing the expedited order $Q^e_t$ given the current state of the system $\Xb^r_t$ as a function of the regular order $Q^r_t$ that we will place in period $t$. 
The ordering rule for the regular supplier can now formally be described as follows. 
In each period $t$ immediately after expedited ordering, we place a regular order $Q^r_t$ such that the expected expedited inventory position $l$ periods later after ordering, $\E[IP^e_{t+l}+Q^e_{t+l} \vert \Xb^r_t]$ equals to the projected expedited inventory position level $\xi$.
That is, the period $t$ regular order $Q^r_t$ solves
\begin{equation} \label{eq:order}
 \E[IP^e_{t+l} (Q^r_t) +Q^e_{t+l} (Q^r_t) \vert \Xb^r_t] = \xi.
\end{equation}
The solution to Equation \eqref{eq:order} is unique and can be efficiently computed using a simple bisection search.  
Proposition~\ref{prop:attainment} in the next subsection shows that there is a non-negative order attaining \eqref{eq:order} in every period after initial attainment.

\subsection{Optimization Procedure}
We now turn to the issue of finding the optimal PEIP policy parameters. Observe that the expedited inventory after ordering is bounded from below by the order-up-to level $S^e$. This implies that the projected expedited inventory position after ordering equals $S^e$ plus the overshoot, i.e., $IP^e_t+Q^e_t=S_e+O_t$. 
A given projected expedited inventory position level $\xi$ thus only scales $IP^e$ for a given $S^e$, which suggests that order placement depends on $S^e$ and $\xi$ only through their difference.
We state this observation formally in the lemma below, which is a special case of Lemma 4.2 of \cite{sheopuri2010newPolicies}:
%\begin{lemma}
%\label{lem:recursions}
%The overshoot $O_t$, the order from the expedited supplier $Q^r_e$, and the order from the regular supplier $Q^r_t$ satisfy the following recursions:
%\begin{align*}
%O_{t+1} &= (O_t+ Q^r_{t-l+1} -D_t)^+,\\
%q_{t+1}^e &= (D_t - O_t -	Q^r_{t-l+1})^+, \\
%Q^r_{t+1} &= \xi - \E[IP^e_{t+1+l} \vert \Xb^r_{t+1}].
%Q^r_{t+1} &= \xi - \E[O_{t+1+l} \vert O_{t+1}, Q^r_{t+2-l}, Q^r_{t+3-l}, \ldots, Q^r_t].
%\end{align*}
%\end{lemma}
%NOTE: deze recursies laten natuurlijk precies zien dat de informatie contained in X^r_t precies de informatie is die nodig is om de overshoot te bepalen. Wellicht moeten we deze recursies naar voren halen. Deze recursies laten ook zien dat voor een l=0, een PEIP precies een PIL is. 

\begin{lemma}
\label{lem:independence}
The distributions of $O$, $Q^r$, and $Q^e$ depend on $\xi$ and $S^e$ only through their difference $V = \xi -S^e$.
\end{lemma}
As our discussion preceding Lemma \ref{lem:independence} already alludes to, the difference $V$ has an explicit meaning that allows for an alternative but equivalent specification and operation of the PEIP policy. 
Rather than keeping the projected expedited inventory position at a fixed level, this alternative specification places regular orders such that the expected overshoot is kept at a fixed level (while retaining an order-up-to level $S^e$ for the expedited
supplier). 
To state this equivalence relation formally, observe that $\E[O_{t+l} \mid \Xb^r_t ]$ denotes the
expected overshoot $l$ periods into the future at time $t$ given the system state at time $t$ following ordering, $\Xb^r_t$.
In any period $t$ after ordering, we have $\E[IP^e_{t+l}+Q^e_{t+l} \mid \Xb^r_t ]=S^e + \E[O_{t+l} \mid \Xb^r_t ]$. 
Hence, placing regular orders such that $\E[IP^e_{t+l}+Q^e_{t+l} \mid \Xb^r_t ]$ equals the projected expedited inventory position level $\xi$ is equivalent to placing regular orders such that $\E[O_{t+l} \mid \Xb^r_t ]$ equals what we henceforth refer to as the projected overshoot level $V=\xi-S^e$. Note also that $O_{t+l}$ is known given $\Xb^r_t$ and $D_t,\ldots,D_{t+l-1}$.

This means in particular that attaining the PEIP level, $\xi$, by placing $Q^r_t$ is equivalent to attaining the projected overshoot level $V$. The next proposition states that attainment is always possible, after the initial attainment.
Note that \cite{vanJaarsveld2024} have established a similar result in the context of lost sales inventory systems.  
The proof of the proposition can be found in Appendix~\ref{sec:attainmentProof}.
\begin{proposition}
\label{prop:attainment}
Under the PEIP policy, if the projected overshoot level $V$ is attainable initially, then in every period the conditional expected overshoot before regular ordering is at most $V$. Consequently, a non-negative regular order attaining $V$ exists in every period.
\end{proposition}
As a convenient notational device let $\xb=(o,q^r_{-l+1},q^r_{-l+2},\ldots q^r_{-1})$ be a generic realization of $\Xb^r_t$ and note $\xb=o$ when $l=1$. For $l\geq 2$, let $\tilde{O}_{0}(\xb)=o$ and
\begin{align*}
    %\tilde{O}_1(\xb)&=(o+q^ r_{-l+1}-D_1)^+,\\
    \tilde{O}_k(\xb)&= \left(\tilde{O}_{k-1}(\xb)+q^r_{-l+k}-D_k \right)^+, \qquad k\in\{1,\ldots,l-1\},\\
    \tilde{O}_l(\xb,y) &= \left(\tilde{O}_{l-1}(\xb)+y-D_l \right)^+.
\end{align*}
Now define the function $g(y \mid \xb)=\E[\tilde{O}_l(\xb,y)]$. Note that $g(y\mid \xb)$ is increasing so it can be inverted. Now define $P_V(\xb)=g^{-1}(V \mid \xb)$ such that $g(P_V(\xb)\mid \xb)=V$. Then $P_V(\Xb^r_t)=Q^r_t$, that is $P_V(\xb)$ returns the regular PEIP order when $\Xb^r_t=\xb$ under projected overshoot level $V$.
We refer to this PEIP policy as $P_V$. 

Observe that, under the PEIP policy, the process $\Xb^r:=\{\Xb^r_t\}_{t\in \N_0}$ is a Markov chain because
$Q^r_t$ is a deterministic function of $\Xb^r_t$. 
Lemma~\ref{lemma:stationarity} (with proof in Appendix~\ref{section:StationarityProof}) establishes that $\Xb^r$ converges weakly to a steady state.
\begin{lemma}
    \label{lemma:stationarity}
    The Markov chain $\Xb^r$ has a compact state-space and a unique invariant probability measure. The invariant measure is attractive so that $\Xb^r_t$ converges in distribution to a stationary state random variable $\Xb^r_\infty$.
\end{lemma}
Let $O_\infty \left(P_V \right)$ and $Q^r_\infty(P_V)$ denote the steady state overshoot and regular order random variables under the PEIP policy $P_V$, respectively.
That is, $O_\infty \left(P_V \right)$ and $Q^r_\infty(P_V)$ are distributed according to the marginals of the invariant probability measure of $\Xb^r$.
The cost-rate of $P_V$, i.e., $C^{P_V}$, can be expressed in terms of $O_\infty \left(P_V \right)$ and $Q^r_\infty(P_V)$ as follows:
\begin{equation*}
\label{eq:SSCostRate_PEIP}
\begin{aligned}
C^{P_V}
    ={}& c^e \left(\mu - \E[Q^r_\infty(P_V)] \right) \\
    &+ h\E \left[ \left(S^e + O_\infty(P_V) - \Db^{l^e+1} \right)^+ \right]
    + p\E \left[ \left(\Db^{l^e+1} - S^e - O_\infty(P_V) \right)^+ \right].
\end{aligned}
\end{equation*}

While the original specification of the PEIP policy has intuitive appeal, the alternative specification in terms of the overshoot outlined above allows us to decompose the optimization of the policy parameters into two simple uni-variate optimization problems: 
\begin{lemma}
\label{lem:optimalS}
The optimal order-up-to level $S^{e*}$ for a given $V$ equals
\[
S^{e*} = \inf \left\{S^e \in \R : \P \left\{ \Db^{l^e+1} - O_\infty \left(P_V \right) \leq S^e \right\} \geq \frac{p}{p+h} \right\}. 
\]
\end{lemma}
The result above implies that a simple search procedure over the projected overshoot level $V$ suffices to find the globally optimal PEIP policy. 

Lemma \ref{lem:optimalS} is a special case of Lemma 4.3 of \cite{sheopuri2010newPolicies}.
\cite{sheopuri2010newPolicies} show that this result holds in general for all policies that have an order-up-to rule for the expedited supplier and operate some given regular ordering rule that uses only the information contained in $\Xb^r_t$. For all such policies, the optimal order-up-to level for the expedited supplier is the solution to a standard Newsvendor model for any policy that places regular orders based only on $\Xb^r$.

\section{Asymptotic Optimality for Long Regular Lead Times}
\label{sec:longleadtime}

In this section, we show that the PEIP policy is asymptotically optimal as the lead time difference $l=l^r-l^e$ approaches infinity, with $l^e$, $h$, $p$, and $c^e$ held fixed. Throughout this section, we assume that $\sigma^2\in(0,\infty)$.

We couple systems operating under the PEIP and TBS policies to show that the cost-rate difference between the two policies decays exponentially as $l$ grows.
The PEIP policy consequently attains an exponentially vanishing optimality gap, as the TBS policy is itself asymptotically optimal at an exponential rate \citep{Chen2026}.

When the lead time difference is $l$, we refer to the PEIP policy with the expedited order-up-to level $S^e$ and $\xi=S^e+V$ as $P_{V,l}$, and to the TBS policy with the same parameter $S^e$ and constant regular order $\varphi\in[0,\mu)$ as $T_\varphi$. We write $C^{P_{V,l}}$ for the cost-rate of the PEIP policy $P_{V,l}$, and $C^{T_\varphi}$ for the cost-rate of the TBS policy $T_\varphi$. Note that costs are independent of $l$ for any TBS policy.

For $\varphi\in[0,\mu)$, the increments $\varphi-D$ are bounded above and have a negative mean.
Therefore, the Lindley recursion $O_{t+1}(T_\varphi)=(O_t(T_\varphi)+\varphi-D_t)^+$ has a stationary distribution \citep[Corollaries III.6.4--III.6.6][]{Asmussen2003}. Observe that this recursion is identical to the waiting-time recursion of a GI/G/1 queue with deterministic service times $\varphi$ and interarrival times distributed as $D_t$. We will exploit this queueing analogy repeatedly in our proofs.
The steady state overshoot random variable under $T_\varphi$, denoted by $O_\infty(T_\varphi)$, is distributed as
\begin{equation*}
    \label{eq:TBS-overshoot}
    O_\infty(T_\varphi)\eqd \sup_{t\in\N_0}\{t\varphi-D[0,t-1]\}.
\end{equation*}
% where $D[0,t-1]$ is cumulative demand over periods $0,\ldots,t-1$ and an empty sum is zero.
Consequently, the cost-rate of $T_\varphi$ is given by
\begin{equation*}
    \begin{aligned}
    C^{T_\varphi}={}&c^e(\mu-\varphi) +h\E \left[\left(S^e+O_\infty(T_\varphi)-\Db^{l^e+1} \right)^+ \right]
    &+p\E \left[ \left(\Db^{l^e+1}-S^e-O_\infty(T_\varphi) \right)^+ \right].
    \end{aligned}
\end{equation*}
% where the demand block $D[0,l^e]$ is independent of the stationary overshoot, as in Equation~\eqref{eq:SSCostRate_PEIP}.
% The supremum representation in \eqref{eq:TBS-overshoot} is a distributional identity; it does not prescribe their joint construction.

\noindent
We restrict attention to $V>0$ and $\varphi > \text{essinf } D_0$, since otherwise $O_\infty(P_{V,l})=O_\infty(T_\varphi)=0$ almost surely.

$O_\infty(T_\varphi)$ is light-tailed \citep[Theorem XIII.5.1][]{Asmussen2003} and all its moments are finite.
The expectation $\E [O_\infty(T_\varphi)]$ is non-decreasing and convex in $ \varphi \in [0,\mu)$, and is strictly increasing where it is positive. Additionally, $\E [O_\infty(T_\varphi)]$ diverges as $\varphi \to \mu$.
Consequently, for every $V>0$, there exists a unique $\varphi_V\in(0,\mu)$ such that
\begin{equation*}
    \label{eq:matching}
    \E[O_\infty(T_{\varphi_V})]=V.
\end{equation*}
Notice that the expected steady state overshoots under $P_{V,l}$ and $T_{\varphi_V}$ are equal, since Proposition~\ref{prop:attainment} gives
\begin{equation*}\label{eq:projected-overshoot-attainment}
 \E[O_{t+l}(P_{V,l})\mid\Xb_t^r]=V,
 \qquad \E[O_\infty(P_{V,l})]=V.
\end{equation*}

Fix $V>0$ and define the generating function $\phi(z):=\E[e^{z(\varphi_V-D)}]$ for $z \geq0$.
\citet{Kingman1962} establishes several properties of the generating function, some of which we use below.
$\phi(z)$ is convex and continuous in $z$.
It is right-differentiable at 0 and the derivative equals to $\varphi_V-\mu<0$.
Since $\varphi_V > \text{essinf } D_0$, we have $\P\{D_0<\varphi_V\}>0$ and consequently $\phi(z)$ diverges as $z \to \infty$. 
Furthermore, $\phi$ attains its minimum at some $z_V>0$ with $\phi(z_V) \in (0,1)$.

Choose $\vartheta \in (0, z_V/2]$ and define
\begin{equation*}\label{eq:convergence-coefficients}
    \gamma:=\phi(\vartheta),\qquad \beta:=\E\!\left[\left(\left(1-\frac{D_0}{\varphi_V}\right)^+\right)^2\right].
\end{equation*}
$\gamma \in (0,1)$ due to convexity of $\phi$ and $\beta>0$ because $\varphi_V> \text{essinf } D_0$. Note also that $\phi(2\vartheta)\leq \phi(\vartheta)=\gamma<1$.
% This cutoff always exists: the transform has right derivative $\phi'(0+)=\varphi_V-\mu<0$ \citep[Proposition XIII.1.1][]{Asmussen2003}, so $\phi(2B)-\phi(B)=(\varphi_V-\mu)B+o(B)<0$ and $\phi(B)<1$ for all sufficiently small $B>0$. Convexity of $\phi$ and $\phi(2B)<\phi(B)$ imply that $\phi$ is strictly decreasing on $[0,B]$. Hence
% \[
%     \vartheta=B>0,\qquad 0<\phi(2\vartheta)<\gamma=\phi(\vartheta)<1.
% \]
% Also, $\beta>0$ because $\varphi_V>\operatorname{ess\,inf}D$.
The coefficients $\vartheta,\gamma,\beta$ depend on the demand distribution and on $V$ through the matching regular order $\varphi_V$.
% $\vartheta$ and $\gamma$ also depend on the chosen cutoff.
We suppress this dependence in the notation.
The resulting coefficients are independent of $l$ and $S^e$.

For the optimal TBS policy $(S^{e^*},\varphi^*)$ where $\varphi^* \in [0,\mu)$ and $S^{e*} \in \R$ let
\[
    V_\infty:=\E[O_\infty(T_{\varphi^*})].
\]
In what follows, we assume that $V_\infty>0$ and omit the trivial case $V_\infty=0$.
Define $\phi^*(z):=\E[e^{z(\varphi^*-D)}]$ with a minimum at $z^* \in (0,\infty)$, and choose $\vartheta^* \in (0,z^*/2]$.
Define
\[
    \gamma^*:=\phi^*(\vartheta^*),\qquad \beta^*:=\E\left[\left(\left(1-\frac{D_0}{\varphi^*}\right)^+\right)^2\right].
\]
Note that $\varphi_{V_\infty}=\varphi^*$ and we use the same generic coefficients at $V=V_\infty$ and the starred coefficients.
% \begin{remark}
%     We select the optimal TBS policy $(S^{e*},\varphi^*)$ as a limit point of optimal TBS pairs for the continuous demand approximations used in Appendix C of \cite{Chen2026}. This pair is optimal by their Lemma EC.7.
% \end{remark}

% Since $\varphi^*<\mu$, the cutoff exists and gives $\vartheta^*>0$ and $0<\phi^*(2\vartheta^*)<\gamma^*<1$, including when $V_\infty=0$.

% All these coefficients are fixed as $l\to\infty$. They need not remain bounded away from their boundary values in a different regime where $\varphi^*\uparrow\mu$.
% The parameter $\vartheta$ determines the exponential weighting of the increments $\varphi_V-D$. The restriction $\phi(2\vartheta)<\gamma<1$ ensures that the stationary TBS overshoot has a finite exponential moment at $2\vartheta$; this allows the conditional-mean difference in Step 2 to have a second-moment bound proportional to $\gamma^{2l}$. Consequently, the mean absolute regular order difference has a bound proportional to $\gamma^l$. The coefficient $\beta$ measures how often, and how far, demand falls below the constant regular order. In those periods the TBS overshoot is not reflected at zero, which allows us to relate changes in the overshoot difference to differences in regular orders. Keeping $\beta$ separately preserves this information in the convergence bounds.
Since $\phi(\vartheta)<1$, $\vartheta$ is smaller than the positive solution of $\phi(z)=1$. Lundberg's inequality \citep[Theorem XIII.5.1][]{Asmussen2003} therefore gives
\begin{equation*}
    \label{eq:TBS-tail}
    \P\{O_\infty(T_{\varphi_V})>x\}\leq e^{-\vartheta x},\qquad x\geq0.
\end{equation*}

Next, we compare two overshoot processes on the same probability space and driven by the same demands and TBS policy.
The first process $\left\{O_{1,t}(T_{\varphi_V}) \right\}_{t \in \N_0}$ is initiated in stationarity.
The second process $\left\{O_{2,t}(T_{\varphi_V}) \right\}_{t \in \N_0}$ starts from $O_{2,0}(T_{\varphi_V})=0$.
The comparison quantifies how quickly a system starting from zero approaches its stationary distribution and how quickly its expected overshoot loses dependence on its initial condition, when the TBS policy is applied.

\begin{lemma}\label{lemma:TBS-mixing}
For every $t\in\N_0$,
\begin{equation}
    \label{eq:TBS-bias}
    0\leq V-\E[O_{2,t}(T_{\varphi_V})]\leq\frac{\gamma^{t+1}}{\vartheta(1-\gamma)},
\end{equation}

\begin{equation}
    \label{eq:TBS-prediction}
    \E\left[\left(\E[O_{1,t}(T_{\varphi_V})\mid O_{1,0}(T_{\varphi_V})]-V\right)^2\right]\leq\frac{\gamma^{2t}}{\vartheta^2(1-\gamma)}.
\end{equation}
\end{lemma}

\proof{Proof of Lemma~\ref{lemma:TBS-mixing}.}
Inequalities \eqref{eq:TBS-bias} follow from the geometric mixing property of the GI/G/1 queue \citep[Theorem 4][]{Kingman1962}.
% The first inequality is the empty-start convergence bound of \cite[Lemma 5][]{xin2018asymptotic}. Its explicit coefficient follows from \cite[Lemma 6][]{Chen2026}. Their Chernoff argument applies to any fixed regular order below $\mu$, any positive parameter with transform below one, and demand distributions with atoms. In our notation, it gives
% \begin{equation}\label{eq:TBS-tail-difference}
%     0\leq\P\{O_\infty(T_{\varphi_V})>x\}-\P\{O_{2,t}(T_{\varphi_V})>x\}\leq\frac{e^{-\vartheta x}\gamma^{t+1}}{1-\gamma},\qquad x\geq0.
% \end{equation}
% Monotonicity of the Lindley recursion gives $O_{1,t}(T_{\varphi_V})\geq O_{2,t}(T_{\varphi_V})$. Write $\Delta O_t^{TBS}:=O_{1,t}(T_{\varphi_V})-O_{2,t}(T_{\varphi_V})$. Since $O_{1,t}(T_{\varphi_V})$ has the stationary distribution, \eqref{eq:TBS-tail-difference} gives
% \begin{align*}
%     \E[\Delta O_t^{TBS}]&=\int_0^\infty\left(\P\{O_{1,t}(T_{\varphi_V})>x\}-\P\{O_{2,t}(T_{\varphi_V})>x\}\right)\,dx\\
%     &\leq\frac{\gamma^{t+1}}{1-\gamma}\int_0^\infty e^{-\vartheta x}\,dx
%     =\frac{\gamma^{t+1}}{\vartheta(1-\gamma)}.
% \end{align*}
% This proves \eqref{eq:TBS-bias}.

Next, we prove Inequality~\eqref{eq:TBS-prediction}.
Define $\Delta O_t^{TBS}:= O_{1,t}(T_{\varphi_V})-  O_{2,t}(T_{\varphi_V})$, for $t \in \N_0$.
Monotonicity of the Lindley recursion implies that
\begin{equation}
    \label{eq:DeltaO_t_nonNegativity}
    \Delta O_t^{TBS} \geq  0,
\end{equation}
since $O_{2,0}(T_{\varphi_V})=0$ and both systems are driven by the same increments $\{\varphi_V-D_n\}_{n \in \N_0}$.
Furthermore, $ \Delta O_t^{TBS}=0$ if and only if in some period $t^\prime \leq t$, $ O_{1,t^\prime}(T_{\varphi_V}) =  O_{2, t^\prime}(T_{\varphi_V})=0$.
This gives
\[
    O_{1,t}(T_{\varphi_V})=\max\{O_{2,t}(T_{\varphi_V}),\,O_{1,0}(T_{\varphi_V})+t\varphi_V-D[0,t-1]\}.
\]
Consequently,
\[
    0\leq\Delta O_t^{TBS}\leq \left(O_{1,0}(T_{\varphi_V})+t\varphi_V-D[0,t-1]\right)^+.
\]
Using $a^+\leq e^{\vartheta a}/\vartheta$ and independence of future demands from the initial stationary overshoot yields
\begin{equation}
    \label{eq:TBS-conditional-difference}
    0\leq\E \left[\Delta O_t^{TBS} \Big \vert O_{1,0}(T_{\varphi_V}) \right] \leq \frac{1}{\vartheta}\E \left[e^{\vartheta\left(O_{1,0}(T_{\varphi_V})+t\varphi_V-D[0,t-1]\right)} \Big \vert O_{1,0}(T_{\varphi_V}) \right]
    \leq\frac{\gamma^t}{\vartheta}e^{\vartheta O_{1,0}(T_{\varphi_V})}.
\end{equation}
Since $O_{2,t}(T_{\varphi_V})$ is independent of $O_{1,0}(T_{\varphi_V})$ and $\E[O_{1,t}(T_{\varphi_V})]=V$,
\begin{alignat}{2}
    &\E\left[\left(\E[O_{1,t}(T_{\varphi_V})\mid O_{1,0}(T_{\varphi_V})]-V\right)^2\right] \nonumber \\
    &= 
    \E\left[\left(\E[O_{1,t}(T_{\varphi_V})-O_{2,t}(T_{\varphi_V})\mid O_{1,0}(T_{\varphi_V})]-\left(V-\E[O_{2,t}(T_{\varphi_V})]\right) \right)^2\right] \nonumber\\
    &=\Var\left(\E[\Delta O_t^{TBS}\mid O_{1,0}(T_{\varphi_V})]\right) \nonumber\\
    &\leq\E\left[\left(\E[\Delta O_t^{TBS}\mid O_{1,0}(T_{\varphi_V})]\right)^2\right] \nonumber\\
    &\leq\frac{\gamma^{2t}}{\vartheta^2}\E[e^{2\vartheta O_{1,0}(T_{\varphi_V})}]. \label{eq:conditional_O2_1}
\end{alignat}
The last inequality in \eqref{eq:conditional_O2_1} follows from \eqref{eq:TBS-conditional-difference}.
The supremum representation \eqref{eq:TBS-overshoot} and $\phi(2\vartheta)\leq \phi(\vartheta)=\gamma<1$ give
\begin{align}
    \label{eq:TBS-exponential-moment}
    \E \left[e^{2\vartheta O_\infty \left(T_{\varphi_V} \right)} \right]&=
     \E \left[e^{2\vartheta \sup_{n\geq 0}\{n\varphi_V-D[0,n-1]\}} \right]
    \leq\sum_{n=0}^{\infty}\E \left[e^{2\vartheta(n\varphi_V-D[0,n-1])} \right] = \sum_{n=0}^{\infty} \phi^n(2 \vartheta) \nonumber\\
    &=\frac{1}{1-\phi(2\vartheta)}\leq\frac{1}{1-\gamma}.
\end{align}
The second equality in \eqref{eq:TBS-exponential-moment} follows from independence and the multiplicative property of moment generating functions under convolution, whereas the third equality and the final inequality follow from $0 < \phi(2\vartheta) \leq \phi(\vartheta) <1$.
Combining \eqref{eq:TBS-conditional-difference} and \eqref{eq:TBS-exponential-moment} yields \eqref{eq:TBS-prediction}.
% For the conditional second moment, we first average over future demands and then square. 
% Here the supremum of the nonnegative exponentials is at most their sum, and independence of the demands evaluates each summand.
% This proves \eqref{eq:TBS-prediction}. By stationarity and time homogeneity, the same inequality holds after any shift of the starting time. Averaging before squaring preserves the factor $\gamma^{2t}$; bounding by $\E[(\Delta O_t^{TBS})^2]$ would lose this improvement.
% \hfill
\Halmos
\endproof

We now couple systems operating under the PEIP policy $P_{V,l}$ and the TBS policy $T_{\varphi_V}$.
They are initialized in their stationary distributions, independently of each other and future demands, and driven with the same demand process.
Lemma~\ref{lemma:stationarity} ensures that this initialization is possible for PEIP. Let
\begin{equation*}
    \label{eq:PEIP-TBS-differences}
    \Delta O_t:=O_t(P_{V,l})-O_t(T_{\varphi_V}),\qquad \DeltaQ_t:=Q_t^r(P_{V,l})-\varphi_V.
\end{equation*}
% PEIP projects its expected overshoot exactly to $V$, whereas the TBS projection is close to $V$ after $l$ periods. Consequently, the expectation $\E[\Delta O_{t+l}\DeltaQ_t]$ is small. We use this specific product in the squared-distance inequality below to control the regular order difference.

\begin{lemma}
    \label{lemma:contraction_identities}
    For $t\geq l$,
    \begin{equation}
        \label{eq:overshoot-contraction}
        \E[(\Delta O_{t+1})^2]\leq \E[(\Delta O_t)^2]-\E[(\Delta O_{t+1}-\Delta O_t)^2] +2\E[\Delta O_{t+1}\DeltaQ_{t-l+1}],
    \end{equation}
    \begin{equation}
        \label{eq:order-identification-bound}
        \E[(\Delta O_{t+1}-\Delta O_t)^2]\geq\beta\E[(\DeltaQ_{t-l+1})^2].
    \end{equation}
\end{lemma}

\proof{Proof of Lemma~\ref{lemma:contraction_identities}}
The overshoot recursions \eqref{eq:Dynamics_O} and $(a^+-b^+)^2\leq(a^+-b^+)(a-b)$ give, for $t\geq l$,
\[
    (\Delta O_{t+1})^2\leq\Delta O_{t+1}(\Delta O_t+\DeltaQ_{t-l+1}).
\]
Multiplying both sides by 2, using $2ab=a^2+b^2-(a-b)^2$ on the right hand side, and taking expectations gives \eqref{eq:overshoot-contraction}.
To derive \eqref{eq:order-identification-bound}, first consider a period with $D_t<\varphi_V$.
In this case $O_t(T_{\varphi_V})+\varphi_V-D_t\geq0$.
If $O_t(P_{V,l})+Q^r_{t-l+1}(P_{V,l})-D_t\geq0$, then $\Delta O_{t+1}-\Delta O_t=\DeltaQ_{t-l+1}$. Otherwise, $Q_{t-l+1}^r(P_{V,l})\leq D_t<\varphi_V$ and
\[
    \vert \Delta O_{t+1}-\Delta O_t\vert =\varphi_V+O_t(P_{V,l})-D_t\geq\varphi_V-D_t,\qquad \vert \DeltaQ_{t-l+1}\vert \leq\varphi_V.
\]
Thus, in either case, and trivially when $D_t\geq\varphi_V$,
\[
    (\Delta O_{t+1}-\Delta O_t)^2\geq\left(\left(1-\frac{D_t}{\varphi_V}\right)^+\right)^2(\DeltaQ_{t-l+1})^2.
\]
Independence of $D_t$ from the order placed at $t-l+1$, including when $l=1$, yields \eqref{eq:order-identification-bound} after taking expectations.
\hfill\Halmos
\endproof

\noindent
It is worth noting that Lemma~\ref{lemma:contraction_identities} does not require the marginal processes to be stationary.

\begin{lemma}\label{lemma:order-convergence}
For every $l\geq1$,
\begin{equation}\label{eq:order-L1}
    \E \left[ \big\vert Q_\infty^r(P_{V,l})-\varphi_V \big \vert \right]\leq\frac{2\gamma^l}{\vartheta\beta\sqrt{1-\gamma}}.
\end{equation}
\end{lemma}

\proof{Proof of Lemma~\ref{lemma:order-convergence}.}
All steady state regular orders have the same marginal distribution, although consecutive orders are not independent. Thus $\E[(\DeltaQ_t)^2]=\E \left[(Q_\infty^r(P_{V,l})-\varphi_V)^2 \right]$ for every $t\geq0$.
Lemma~\ref{lemma:stationarity} and Inequality~\eqref{eq:TBS-tail} ensure that all second moments below are finite for each fixed $l$.

Future demands are independent of the two current states. Conditioning on these states and using \eqref{eq:projected-overshoot-attainment} gives
\[
    \E[\Delta O_{t+l}\mid\Xb_t^r,O_t(T_{\varphi_V})]=V-\E[O_{t+l}(T_{\varphi_V})\mid O_t(T_{\varphi_V})].
\]
Because $\DeltaQ_t$ is determined by $\Xb_t^r$, Cauchy--Schwarz and \eqref{eq:TBS-prediction} imply
\begin{alignat}{2}
    \label{eq:projection-product}
    & \vert \E[\Delta O_{t+l}\DeltaQ_t] \vert = \vert \E\left[  \E \left[\Delta O_{t+l} \mid \Xb_t^r,O_t(T_{\varphi_V})\right] \DeltaQ_t\right] \vert \nonumber\\
    & \leq \left(\E \left[ \left( \E \left[\Delta O_{t+l} \mid \Xb_t^r,O_t(T_{\varphi_V})\right]\right)^2 \right] \right)^{1/2} \left(\E[(\DeltaQ_t)^2]\right)^{1/2}
    \leq \frac{\gamma^l}{\vartheta\sqrt{1-\gamma}}\left(\E[(\DeltaQ_t)^2]\right)^{1/2}.
\end{alignat}

% The last expectation is bounded by \eqref{eq:projection-product} applied at time $t-l+1$.

\noindent
Combining \eqref{eq:overshoot-contraction}, \eqref{eq:order-identification-bound}, and \eqref{eq:projection-product}, summing over $t=l,\ldots,l+N-1$, and dividing by $N$, we obtain by stationarity of the regular orders
\[
    \beta\E \left[(Q_\infty^r(P_{V,l})-\varphi_V)^2 \right]  \leq\frac{\E[(\Delta O_l)^2]}{N}+\frac{2\gamma^l}{\vartheta\sqrt{1-\gamma}}\left(\E \left[(Q_\infty^r(P_{V,l})-\varphi_V)^2 \right]\right)^{1/2}.
\]
Here the nonnegative terminal squared distance has been dropped. Let $N\to\infty$, keeping $l$ fixed.
Then
\begin{equation*}
    \label{eq:order-L2}
    \E[(Q_\infty^r(P_{V,l})-\varphi_V)^2]\leq\frac{4\gamma^{2l}}{\vartheta^2\beta^2(1-\gamma)}.
\end{equation*}
% If $\E[(\DeltaQ_0)^2]=0$, the result is immediate; otherwise divide by its square root and square to obtain \eqref{eq:order-L2}.
Since $f(x)=x^2$ is a convex function, Jensen's inequality yields
\[
    \left(\E \left[\vert Q_\infty^r(P_{V,l})-\varphi_V \vert \right]\right)^2 \leq \E[(Q_\infty^r(P_{V,l})-\varphi_V)^2]\leq\frac{4\gamma^{2l}}{\vartheta^2\beta^2(1-\gamma)}
\]
Taking square roots then gives \eqref{eq:order-L1}.
\hfill\Halmos
\endproof

Lemma~\ref{lemma:order-convergence} can be extended to the cost-rates.
The key step in this extension is to prove 
that $O_\infty(P_{V,l})$ converges in distribution to $O_\infty(T_{\varphi_V})$ as $l$ grows
This is required because the cost-rates depend on the full distributions of the corresponding overshoot random variables.

\begin{theorem}\label{theorem:cost-convergence}
$O_\infty(P_{V,l})$ converges in distribution to $O_\infty(T_{\varphi_V})$ as $l \to \infty$.
Moreover for $l\geq1$,
\begin{equation}
    \label{eq:cost-convergence-bound}
     \Big\vert C^{P_{V,l}}-C^{T_{\varphi_V}} \Big\vert \leq \frac{2 \gamma^l}{\vartheta \sqrt{1-\gamma}} \left(\frac{c^e + 2l\max\{h,p\}}{\beta} + \frac{\gamma \max\{h,p\}}{\sqrt{1-\gamma}} \right)
\end{equation}
uniformly over the common expedited order-up-to level $S^e\in\R$. 
\end{theorem}

\proof{Proof of Theorem~\ref{theorem:cost-convergence}.}
We couple two systems on a common probability space. The first operates under the PEIP policy and is initialized in stationarity, whereas the second operates under the TBS policy with $O_l(T_{\varphi_V})=0$ at time $t=l$.
Both systems are driven by the same demand sequence.
$\Delta O_t$ and $\DeltaQ_t$ are defined according to \eqref{eq:PEIP-TBS-differences}.
Using $\vert a\vert =a+2(-a)^+$ we have
\begin{equation}
    \label{eq:absoluteDelta_O}
     \E[\vert \Delta O_t \vert] =  \E[\Delta O_t] + 2\E[(-\Delta O_t)^+].
\end{equation}
In the following we bound the two terms on the right hand side of \eqref{eq:absoluteDelta_O}.
By \eqref{eq:TBS-bias}, shifted to start at time $l$, we have
\begin{equation}
    \label{eq:absoluteDelta_O_1}
    0\leq\E[\Delta O_t]\leq\frac{\gamma^{t-l+1}}{\vartheta(1-\gamma)}.
\end{equation}
Next, observe that, for $t\geq l$, the positive-part inequality gives
\[
    (-\Delta O_{t+1})^+\leq(-\Delta O_t)^++(-\DeltaQ_{t-l+1})^+.
\]
As noted, $O_l(T_{\varphi_V})=0$ so that $(-\Delta O_l)^+=0$. Now after iterating and summing, stationarity of the PEIP orders, and taking expectations, we obtain
\begin{equation}
    \label{eq:finite-window-overshoot}
    \E[(-\Delta O_t)^+]\leq(t-l)\E \left[ \big\vert Q_\infty^r(P_{V,l})-\varphi_V \big \vert \right].
\end{equation}
It follows from \eqref{eq:order-L1}, \eqref{eq:absoluteDelta_O}, \eqref{eq:absoluteDelta_O_1}, and \eqref{eq:finite-window-overshoot} that
\begin{equation}
    \label{eq:L1_O}
    \E[\vert \Delta O_t \vert]\leq\frac{\gamma^{t-l+1}}{\vartheta(1-\gamma)}+ \frac{4(t-l)\gamma^l}{\vartheta\beta\sqrt{1-\gamma}}.
\end{equation}
For any $1-$Lipschitz function $f$, Inequality~\eqref{eq:DeltaO_t_nonNegativity} together with \eqref{eq:TBS-bias} and \eqref{eq:L1_O} give
\[
    \vert \E[f(O_t(T_{\varphi_V}))]-\E[f(O_\infty(T_{\varphi_V}))]\vert \leq 
    \vert \E[O_t(T_{\varphi_V})]-\E[O_\infty(T_{\varphi_V})] \vert
    \leq\frac{\gamma^{t-l+1}}{\vartheta(1-\gamma)},
\]
and
\[
    \vert \E[f(O_t(T_{\varphi_V}))]-\E[f(O_\infty(P_{V,l}))]\vert \leq\frac{\gamma^{t-l+1}}{\vartheta(1-\gamma)}+ \frac{4(t-l)\gamma^l}{\vartheta\beta\sqrt{1-\gamma}}.
\]
Combining these two inequalities and taking $t=2l$ gives
\begin{equation}
    \label{eq:stationary-overshoot-comparison}
    \vert\E[f(O_\infty(P_{V,l}))]-\E[f(O_\infty(T_{\varphi_V}))]\vert \leq\frac{2\gamma^{l+1}}{\vartheta(1-\gamma)}+ \frac{4l\gamma^l}{\vartheta\beta\sqrt{1-\gamma}}.
\end{equation}
The right hand side of \eqref{eq:stationary-overshoot-comparison} goes to zero as $l$ approaches infinity.
Applying this inequality to bounded Lipschitz functions proves that $O_\infty(P_{V,l})$ convergence in distribution to $O_\infty(T_{\varphi_V})$ by the Portmanteau theorem.
The expected holding and backlog cost, as a function of the overshoot, is Lipschitz with constant $\max\{h,p\}$. The stationary cost formulas, $\vert\E[Q_\infty^r(P_{V,l})]-\varphi_V\vert\leq\E[\vert Q_\infty^r(P_{V,l})-\varphi_V\vert]$, and \eqref{eq:stationary-overshoot-comparison} imply
\[
    \vert C^{P_{V,l}}-C^{T_{\varphi_V}}\vert\leq c^e \E[\vert Q_\infty^r(P_{V,l})-\varphi_V\vert] + \max\{h,p\} \left( \frac{2\gamma^{l+1}}{\vartheta(1-\gamma)}+ \frac{4l\gamma^l}{\vartheta\beta\sqrt{1-\gamma}}\right).
\]
This result combined with \eqref{eq:order-L1} gives \eqref{eq:cost-convergence-bound}.
\Halmos
\endproof

\noindent
The leading term in the bound in \eqref{eq:cost-convergence-bound} is of order $l\gamma^l$, which decays exponentially because $\gamma<1$.
Moreover, it is independent of $S^e$, that is it also bounds the absolute difference between the two cost-rates after it is minimized over $S^e$.

We conclude this section by establishing asymptotic optimality of the PEIP policy with an exponential decay rate in $l$. Define $C^{P_{V^*}}(l):=\inf_{V>0,S^e\in\R}C^{P_{V,l}}$ and $C^{TBS^*}:=\inf_{0\leq\varphi<\mu,S^e\in\R}C^{T_\varphi}$, and let $C^{u^*}(l)$ be the optimal cost-rate over all admissible policies when the lead time difference is $l$.
We use the optimal TBS parameters $(S^{e*},\varphi^*)$, the matching overshoot $V_\infty$, and the starred coefficients defined after \eqref{eq:convergence-coefficients}.

\begin{theorem}\label{theorem:asymptotic-optimality}
% If $V_\infty>0$, then, 
For every $l\geq1$,
\begin{alignat}{2}
    \label{eq:absolute-optimality-bound}
    C^{P_{V^*}}(l)-C^{u^*}(l)\leq& \frac{2 \left(\gamma^* \right)^l}{\vartheta^* \sqrt{1-\gamma^*}} \left(\frac{c^e + 2l\max\{h,p\}}{\beta^* } + \frac{\gamma^* \max\{h,p\}}{\sqrt{1-\gamma^*}} \right)\nonumber\\
    &+\frac{(h+p)\mu\,e^{\vartheta^*[S^{e*}-(l^e+1)\varphi^*]}}{(1-\gamma^*)^2}\,(\gamma^*)^{l^r+2}.
\end{alignat}
% If $V_\infty=0$, then, for every $l\geq1$,
% \begin{equation}\label{eq:absolute-optimality-bound-zero}
%     0\leq C^{P_{V^*}}(l)-C^{u^*}(l)\leq\frac{(h+p)\mu\,e^{\vartheta^*[S^{e*}-(l^e+1)\varphi^*]}}{(1-\gamma^*)^2}\,(\gamma^*)^{l^r+2}.
% \end{equation}
% Consequently, with demand, $l^e$, $h$, $p$, and $c^e$ held fixed, the absolute optimality gap is $O(l(\gamma^*)^l)$ and converges to zero exponentially fast as $l\to\infty$.
\end{theorem}

\proof{Proof of Theorem~\ref{theorem:asymptotic-optimality}.}
Choose an optimal TBS pair $(S^{e*},\varphi^*)$ as provided by Appendix C of \cite{Chen2026}.
Theorem 2, Corollary 2, and the final absolute bound in the proof of Theorem EC.1 of \cite{Chen2026} give
\begin{equation}\label{eq:TBS-optimality-gap}
    0\leq C^{TBS^*}-C^{u^*}(l)\leq\frac{(h+p)\mu\,e^{\vartheta^*[S^{e*}-(l^e+1)\varphi^*]}}{(1-\gamma^*)^2}\,(\gamma^*)^{l^r+2}.
\end{equation}
% The smaller fixed interval $[0,B]$ used here is also admissible in that argument: their Lemma 6 uses only $\vartheta^*>0$ and $\phi^*(\vartheta^*)=\gamma^*<1$, and the convergence results in Equations (EC.15)--(EC.16) are uniform on every fixed compact interval.
% Thus, the selected optimal pair and the coefficients defined above satisfy \eqref{eq:TBS-optimality-gap}.
Use the PEIP policy with projected overshoot $V_\infty$ and expedited order-up-to level $S^{e*}$.
Then
\[
    0\leq C^{P_{V^*}}(l)-C^{u^*}(l)\leq\vert C^{P_{V_\infty,l}}-C^{T_{\varphi^*}}\vert+C^{TBS^*}-C^{u^*}(l).
\]
Theorem~\ref{theorem:cost-convergence} and \eqref{eq:TBS-optimality-gap} give \eqref{eq:absolute-optimality-bound}.
\hfill\Halmos
\endproof

We conclude this section by noting that our proof technique also
applies to the PIL policy for lost sales inventory systems
\citep{vanJaarsveld2024}. As discussed in Section~\ref{sec:literatureReview}, \cite{sheopuri2010newPolicies} show that a lost sales inventory system can be viewed as a special case of a dual-sourcing system when the expedited supplier follows an order-up-to rule. In particular, the pre-receipt inventory in a lost sales
system follows the same recursion as the overshoot in
\eqref{eq:Dynamics_O}.
This allows a similar coupling argument, with a constant order
policy replacing the TBS benchmark. The PIL target is then set
equal to the expected steady state inventory immediately after
receipt under the constant order policy.

After establishing the corresponding steady state result for PIL akin to the stationarity of PEIP,
we may compare the stationary PIL system over one lead time with a
constant order system initialized in stationarity independently
of the PIL history. Driving both systems with the same subsequent
demands gives a contraction inequality analogous to
\eqref{eq:overshoot-contraction}. In this case, the decrease in
the expected squared inventory difference directly bounds the
expected squared order difference. Summing over the comparison
interval and using a lead-time-uniform second-moment bound on
the PIL inventory shows that its orders converge to the constant
order. Target attainment and the steady state inventory balance
then give convergence of the cost-rates. Choosing an optimal
constant order benchmark therefore establishes asymptotic
optimality of PIL in the lost sales system for long leadtimes \citep{xin2016asymptotic}. More generally, our
proof technique may serve as a template for other projection
policies for which similar steady state, coupling, and
contraction inequalities can be established.

% Theorem~\ref{theorem:asymptotic-optimality} separates the difference $C^{TBS^*}-C^{u^*}(l)$ from the additional difference $\vert C^{P_{V_\infty,l}}-C^{T_{\varphi^*}}\vert $ when $V_\infty>0$. With the same admissible coefficients, the bound of \cite[Theorem 2 and Theorem EC.1][]{Chen2026} is proportional to $(\gamma^*)^{l^r+2}$, while our additional PEIP comparison has leading order $l(\gamma^*)^l$. For fixed $l^e$, the two bounds therefore have the same exponential base, with an additional factor $l$ in the PEIP comparison from accumulating regular order differences over an interval of length $l$. The coefficient $\beta^*$ enters the prefactor through the PEIP contraction argument and has no counterpart in their TBS comparison. When $V_\infty=0$, \eqref{eq:absolute-optimality-bound-zero} gives the same bound as the TBS optimality bound.

% The restriction on the cutoff is needed for our conditional second-moment argument. Choosing a larger cutoff in the TBS bound can produce a smaller value of $\gamma^*$, so the comparison above uses the same coefficients in both bounds and does not assert equality with the fastest TBS rate obtainable from their argument. In particular, the revised $\gamma^*$ need not equal a coefficient computed before imposing $\phi^*(2\vartheta^*)<\gamma^*$. The bound $O(l(\gamma^*)^l)$ includes a polynomial prefactor; it does not assert a bound by a constant times $(\gamma^*)^l$. These comparisons concern the proved upper bounds and do not establish an ordering of the actual policy gaps.

\section{Asymptotic Optimality in the High Service Regime}
\label{sec:asymptoticOptimality}
In this section, we prove that the best PEIP policy is asymptotically optimal when  both  the  shortage cost  and  the  cost  premium  for  expedited  units  become  large,  with  their  ratio  held constant.
We first elaborate on why it is important and non-trivial to study this asymptotic regime. We then present the asymptotic optimality result itself, and we subsequently show that this result is robust in the sense that it holds true for a broad class of policies. %generalize this result

\subsection{The High Service Asymptotic Regime}
We know from \cite{xin2018asymptotic} that a simple TBS policy is asymptotically optimal as the lead time
difference between the regular supplier and the expedited supplier approaches infinity. Unfortunately, the policy
performs poorly outside this asymptotic regime, particularly when either the cost for a unit in backlog or the price charged by the expedited supplier is large. 
When these costs grow large simultaneously, the performance of the TBS drops precipitously \citep{klosterhalfen2011comparisonDIP}. 
It is however exactly this latter asymptotic regime that is prevalent in practice, and thus establishing asymptotic performance guarantees of the PEIP policy in this regime has important practical implications.
In the capital goods industry, for instance, backordering a critical part comes at a high cost while expediting shipments from an emergency supplier are expensive too \citep[see, e.g.,][]{westerweel2018printing}.  

In this section, we thus focus on the asymptotic regime where the unit expediting cost $c^e$ and the shortage cost $p$ grow large simultaneously with their ratio held fixed at a constant strictly smaller than $l$ (i.e. $c^e< p\cdot l$). We call this the high service regime.
Recall from Section \ref{sec:description} that if we would not impose this assumption, then the dual-sourcing inventory problem would reduce to a trivial single supplier inventory problem. 

We remark that, following a reasoning similar to the above, the asymptotic regime where both $l$ and $c^e$ grow large with their ratio held fixed at a constant strictly smaller than $p$ is, at least theoretically, interesting as well.
Indeed, this too is a non-trivial regime in which the optimal policy does not reduce to a trivial single supplier inventory problem. \cite{xin2020approx} has recently made advances on this latter regime by showing that for the continuous time variant of the problem, the CDI policy has a worst-case performance guarantee of 1.79.

\subsection{Asymptotic Optimality}
To facilitate our analysis, we first introduce additional notation. Let $C^{P_{V^*}}(p,c^e,l^r,l^e)$ and $C^{u^*}(p,c^e,l^r,l^e)$ denote the cost-rate of the dual-sourcing inventory system under the best PEIP policy and the optimal policy, respectively, as functions of the problem specific input parameters.
Here, the best PEIP policy is taken over $\xi\in\R$ and $S^e\in \R\cup\{-\infty\}$, where $S^e=-\infty$ means that no expedited orders are placed.
We may occasionally omit certain input parameters when there is no ambiguity. Throughout this section, the holding cost $h$ and both lead times $l^r$ and $l^e$ are held fixed. We define $p(n):=p\cdot n$ and $c^e(n):=c^e\cdot n$. 
%For later use we also define $\bar c^e:=\min\{c^e,p/2\}$ and $\bar c^e(n):=\min\{c^e(n),p(n)/2\}=\bar c^e \cdot n.$
We are interested in the asymptotic regime where $n\to\infty$ while
$c^e(n)/p(n)=c^e/p$ remains fixed at a constant strictly smaller than $l$. For a given $n>0$, we further define
$C^{P_{V^*}}(n) = C^{P_{V^*}}(p(n),c^e(n))$ and $C^{u^*}(n) = C^{u^*}(p(n),c^e(n))$ to denote the cost as we scale $p$ and $c^e$ proportionally with $n$. 

Let $\mathcal{L}(p,l)$ and $\mathcal{B}(p,l)$ denote the canonical  single supplier inventory systems with lost sales and backlogged demand, respectively, as defined in
\cite{janakiraman2007comparison}, \cite{huh2009asymptotic}, and \cite{bijvank2014robust}. 
Following their notation, let $C^{\mathcal{L}^*}(p,l)$ denote the cost-rate under an optimal policy in $\mathcal{L}(p,l)$ as a function of the lost sales penalty cost $p$ and the lead time $l$.
Similarly, let $C^{\mathcal{B}^*}(p,l)$ denote the cost-rate under the optimal policy in $\mathcal{B}(p,l)$, with $p$ now denoting the penalty cost for backlogged demand. 
It is well-known that the optimal policy for $\mathcal{B}(p,l)$ is a simple order-up-to rule. It is equally well-known that the optimal policy for $\mathcal{L}(p,l)$ has a complex structure and depends in general on the entire vector of outstanding orders, very much akin to the present dual-sourcing inventory system.

We first establish a lower bound on the cost-rate of the optimal policy, i.e. $C^{u^*}(p,c^e,l^r,l^e)$, in terms of both the cost-rate of the optimal policy for $\mathcal{B}(p,l)$, i.e. $C^{\mathcal{B}^*}(p,l)$, and the optimal policy for $\mathcal{L}(p,l)$, i.e. $C^{\mathcal{L}^*}(p,l)$.  
\begin{lemma}
\label{lem:lowerbound} 
For any $\varepsilon\in(0,1)$, 
\[C^{\mathcal{B}^*}(\min(c^e,\varepsilon p)/(l^r+1),l^r)
\leq
C^{\mathcal{L}^*}(\min(c^e,\varepsilon p)/(l^r+1),l^r)
\leq
C^{u^*}(p,c^e,l^r,l^e)\]
\end{lemma}

\proof{Proof.}
Since the cost-rate of every fixed admissible policy is
non-decreasing in $c^e$, the optimal cost-rate is also
non-decreasing in $c^e$. For some $\varepsilon\in(0,1)$, relpace $c^e$ with $\bar{c}^e=\min\{c^e,\varepsilon p\}<p$. This can only decrease $C^{u^*}$, so it suffices to establish a lower bound for the
system with expedited unit price $\bar c^e$. Observe that a lower bound on the cost-rate of the optimal policy is given by the cost-rate cost of an optimal policy for a system identical to the present inventory system but in which any backlogged demand at the end of a period is satisfied instantly by
the expedited supplier at unit price $\bar c^e$. According to Theorem 3.1 of
\cite{sheopuri2010newPolicies}, this system is exactly
$\mathcal{L}(\bar{c}^e,l^r)$, where the lost sales penalty cost
equals $\bar{c}^e$.
Hence, we have
$C^{u^*}(p,c^e,l^r,l^e)
\geq C^{\mathcal{L}^*}(\bar c^e,l^r)$.
We further know from Theorem 5 of
\citet{janakiraman2007comparison} that
$C^{\mathcal{L}^*}(\bar c^e,l^r)
\geq C^{\mathcal{B}^*}(\bar c^e/(l^r+1),l^r)$,
and thus
$C^{u^*}(p,c^e,l^r,l^e)
\geq C^{\mathcal{B}^*}(\bar c^e/(l^r+1),l^r)$.
\hfill\Halmos
\endproof

We note in passing that easily computable lower bounds for the optimal cost of dual-sourcing systems are relatively scarce. In fact, we are aware of only one such lower bound in terms of the optimal policy for an equivalent dual-sourcing system where the lead time difference is 1 \citep{sheopuri2010newPolicies}. For this equivalent system, the optimal policy is a dual-index policy.
In this respect, the lower bound in terms of $\mathcal{B}(\min(c^e,\varepsilon p)/(l^r+1),l^r)$ in Lemma \ref{lem:lowerbound} is of independent interest as it is easy to compute and leads to an upper bound on the optimality gap for any given policy that increases with $\varepsilon$. 
%We briefly reflect on the quality of this lower bound in our numerical investigation in the next section. 

The result in the high service regime requires the following technical assumption regarding the cumulative demand over $l^r+1$ periods:
\begin{assumption}
\label{demandAssumption}
The random variable $\Db^{l_r+1}$ satisfies one of the following conditions: 
\begin{enumerate}[label=\roman*.]
\item $\Db^{l_r+1}$ is bounded, or
\item $\Db^{l_r+1}$ is unbounded and $\lim_{d\to\infty} \E[\Db^{l_r+1}-d \mid \Db^{l_r+1}>d ] /d =0$. 
\end{enumerate}
\end{assumption}
This assumption first appeared in \cite{huh2009asymptotic}. While this assumption is stronger than the finite-second-moment assumption imposed on demand in the long lead time regime of Section \ref{sec:longleadtime}, \cite{huh2009asymptotic} show that many commonly used distributions in
inventory models satisfy this assumption, including Poisson distributions, geometric distributions, negative
binomial distributions with parameters $r>0$ and $0<p<1$, and (log-)normal distributions. A sufficient condition for
this assumption to be satisfied is that the demand distribution has an increasing failure rate. 

We can now state the main result of this section.
\begin{theorem}
\label{asymptoticallyOptimal}
Let the backorder penalty cost $p$ and the expedited cost premium $c^e$ be such that $c^e/p< l$. Then, under Assumption \ref{demandAssumption}, the best PEIP policy is asymptotically optimal as $c^e$ and $p$ approach infinity simultaneously, that is
\[
\lim\nolimits_{n \to \infty} \frac{C^{P_{V^*}}(n)}{C^{u^*}(n)} = 1.
\]
\end{theorem}

\proof{Proof of Theorem~\ref{asymptoticallyOptimal}.}
Observe that single sourcing from the regular supplier is a feasible PEIP policy by setting $S^e$ equal to $-\infty$.
This observation implies that the numerator is
bounded from above by the cost-rate of an optimal policy in a single supplier inventory system
with backlogged demand whose lead time is equal to the lead time of the regular supplier, i.e. $C^{\mathcal{B}^*}(p,l^r)\geq C^{P_{V^*}}(p,c^e,l^r,l^e)$. Combining this observation with Lemma \ref{lem:lowerbound}, we have using $\bar{c^e}:=\min(c^e,\varepsilon p)$ for some $\varepsilon\in(0,1)$
\[
\frac{C^{P_{V^*}}(n)}{C^{u^*}(n)}
\leq
\frac{C^{\mathcal{B}^*}(p(n),l^r)}
{C^{\mathcal{B}^*}(\bar c^e(n)/(l^r+1),l^r)}.
\]
By Theorem 2b of \cite{huh2009asymptotic}, we know that under Assumption \ref{demandAssumption}
\[
\lim_{n\to\infty}
\frac{C^{\mathcal{B}^*}(p(n),l^r)}
{C^{\mathcal{B}^*}(\bar c^e(n)/(l^r+1),l^r)}
=1.
\]
which gives the desired result.   
\hfill \Halmos \endproof

Define the subset $\tilde{\UC} = \{ u \in \UC \vert  C^{u}(p,c^e,l^r,l^e) \leq C^{\mathcal{B}^*}(p,l^r) \}$. Note that it readily follows from the definition of $\tilde{\UC}$ and the proof of Theorem \ref{asymptoticallyOptimal}, that any policy $u\in\tilde{\UC}$ is asymptotically optimal in the regime under our consideration. 
This is an important observation that generalizes our result to a wide class of policies. We therefore state it as a theorem. 
\begin{theorem}
\label{robust}
Let the backorder penalty cost $p$ and the expedited cost premium $c^e$ be such that $c^e/p< l$. Then, under Assumption
\ref{demandAssumption}, any policy $\tilde{u}\in\tilde{\UC}$ is asymptotically optimal as $c^e$ and $p$ approach infinity simultaneously, that is
\[
 \lim\nolimits_{n \to \infty} \frac{C^{\tilde{u}^*}(n)}{C^{u^*}(n)} =  1,
\]
where $C^{\tilde{u}^*}(n)$ is defined in a way similar as $C^{u^*}(n)$.  
\end{theorem}

We conclude this section by noting that the CDI, DI, and SI policy are all contained in the set $\tilde{\UC}$. They are thus asymptotically optimal when both $c^e$ and $p$ grow large with their ratio held fixed at a constant strictly smaller than $l$.

\section{Numerical Investigation}
\label{sec:numericalStudy}
The theoretical results of the previous sections are important as they provide performance guarantees for the PEIP policy in long lead time and high service regimes. 
While this is particularly useful for practitioners operating in those regimes, it does not imply that the PEIP policy consistently performs well, especially outside those asymptotic regimes.    
In this section, we therefore undertake an extensive numerical investigation to examine the overall effectiveness of the PEIP policy. This investigation mainly revolves around the comparison of the performance of the PEIP policy with the performance of several existing heuristic policies.

We first describe the benchmark heuristic policies in Section \ref{subsec:benchmarks}. Then Section \ref{subsec:comparison}, compares the performance of the PEIP policy with the performance of the benchmark heuristic policies based on a full factorial test bed of industrial size. 

\subsection{Benchmark Heuristic Policies and Computational Aspects}
\label{subsec:benchmarks}
The following five heuristic policies will serve as benchmarks: (i) The  DI policy, (ii) the CDI policy, (iii) the VBS policy, (iv) the BWB policy, and (v) the TBS policy. 
Although we have discussed these five heuristic policies in our literature review, we revisit a few points here that are important to bear in mind while considering the results of the numerical investigation. %discuss how they relate to implemention

First, the CDI policy and the BWB policy are the current best performing heuristic policies. Their computational requirements are however significantly higher than the PEIP policy.
Indeed, while the PEIP policy requires only a one-dimensional search to find the optimal policy parameters, both the CDI policy and the BWB policy require a search over two dimensions to arrive at their best policy parameters.
By contrast, the DI policy, VBS policy, and the TBS policy have the same computational requirements as the PEIP policy.
Second, the CDI policy is a generalization of the DI policy, and it is therefore to be expected that the former policy outperforms the latter. The same holds for the BWB policy and the VBS policy -- the latter being a special case of the former. %maybe add here asymptotic again and robust optimization optimal CDI
Third, all five benchmark heuristic policies use an order-up-to rule for the expedited supplier. 
As such, echoing our discussion after Lemma \ref{lem:optimalS}, the optimal order-up-to level under each policy for given regular ordering rule parameters follows a Newsvendor type equation. 
This is convenient numerically as it only requires us to simulate the steady state overshoot distribution for given regular ordering rule parameters to arrive at the corresponding optimal order-up-to level. 
We describe this simulation-based optimization scheme in detail next.  

For all problem instances in this numerical investigation, we find the best policy parameters for the PEIP policy as well as for the benchmark heuristic policies via golden-section searches over their respective regular ordering rule parameters. 
In each iteration of this search procedure, we first compute the steady state overshoot distribution via simulation, we then determine the optimal order-up-to level via the Newsvendor equation, and we finally compute the corresponding cost-rate. 
We use common random numbers across all policies per problem instance, and we continue simulating the steady state overshoot distribution until the width of the 95 percent
confidence intervals for the mean and for the standard deviation of the overshoot are less than 1 percent
of their respective point estimates. 
After we find the best control policy parameters for all heuristic policies for a given problem instance, we compute the cost-rate for each heuristic policy under its best policy parameters via simulation.
As in the simulation-based optimization scheme, this simulation uses common random numbers across policies, and it continues until the width of the 95 percent confidence intervals for the cost-rate is less than 1 percent
of its point estimate.

\subsection{Full Factorial Experiment}
We compare the performance of the PEIP policy with the performance of the heuristic policies described in the previous section based on a large test bed of industrial size. 
Note that the optimal policy cannot be computed for such large instances due to the complexity of the corresponding dynamic program. 
We therefore exclude the optimal policy in this investigation. 

In all instances of the test bed, the one period demand has a negative binomial distribution with mean 50 and the holding cost $h$ is fixed at 1. We further keep the price of the regular supplier $c^r$ fixed at 0 since only the cost difference between the suppliers affects the relative cost-rate. 
The expedited lead time $l^e$ is kept fixed at 0 periods so that the focus is on the lead time difference $l=l^r-l^e=l^r$. 
The remaining input parameters are varied over multiple levels. We vary the backorder penalty cost over five levels ($p\in\{4,9,19,49,99\}$), the coefficient of variation of the one period demand over six levels ($\sigma/\mu\in\{0.15,0.25,0.5,1,1.5,2\}$), and the regular lead time over four levels ($l^r\in\{2,3,4,5\}$). We set the price of the expedited supplier $c^e$ equal to a fraction of the backorder penalty cost multiplied by the lead time difference of the specific instance, that is, $c^e=\delta pl$ for $\delta\in\{0.1,0.2,0.4\}$.  
While this approach is different from most numerical studies on dual-sourcing inventory systems, it does ensure that there is no instance in the test bed for which exclusively sourcing from the regular supplier is optimal.  
Table \ref{tab:testbed1} summarizes all levels for each input parameter in our test bed that is not kept fixed. Permuting over all these levels leads to a large test bed of 360 instances. 

\label{subsec:comparison}
\begin{table*}[!htbp]
  \centering
  \caption{Input parameter values for our test bed}
	\fontsize{8pt}{9pt}\selectfont
  \label{tab:testbed1}
  \begin{tabularx}{0.71\textwidth}{c X c l}
    \toprule
    & Input parameter& No. of choices&Values\\
		\midrule
		1 & Backorder penalty cost, $p$ & 5 & 4, 9, 19, 49, 99\\
		2 & Regular lead time, $l^r$ & 4 & 2, 3, 4, 5 \\
		3 & Coefficient of variation of demand, $\sigma/\mu$ & 6 & 0.15, 0.25, 0.5, 1, 1.5, 2\\
		4 & Price of the expedited supplier, $c^e$ & 3 & $\delta  p l$ for $\delta$ = 0.1, 0.2, 0.4 \\
    \bottomrule
  \end{tabularx}
\end{table*}

We are interested in the performance of the PEIP policy in comparison with the performance of the benchmark heuristics. 
To that end, we compute for each instance of the test bed the percentage gap in the cost-rate of each benchmark heuristic policy $u\in\{\text{DI,CDI,VBS,BWB,TBS}\}$ relative to the PEIP policy:  $100\cdot \frac{C^{u^*} - C^{P_{V^*}}}{C^{P_{V^*}}}$.

Table \ref{tab:resultsOverall} summarizes the numerical results of the full factorial experiment.
The table presents the average, maximum, and minimum gap percentages for each benchmark heuristic policy computed over all instances as well as the percentage of test instances in
which the PEIP policy outperforms the specific benchmark heuristic policy. (A positive gap percentage implies that the PEIP policy outperforms the specific benchmark heuristic policy.) Detailed results are relegated to Appendix \ref{sec:detailedResults}.     
Based on this table, we conclude that the PEIP policy outperforms the TBS, DI, CDI, and VBS policy. The performance of the PEIP policy seems to be comparable with the BWB policy, though on average the PEIP policy outperforms this policy as well. It is also interesting to note that the PEIP policy can outperform the BWB and CDI policy by up to 2.67\% and 4.67\%, respectively. 
The opposite never holds. Indeed, when the PEIP policy is outperformed by the BWB or CDI policy, it only occurs by much smaller margins of up to 0.29\% and 0.08\%, respectively. Note that these gaps are smaller than
the confidence interval around the estimated performance for any single instance.
This is an interesting observation, especially in view of the fact that the BWB and CDI policy are the current best performing heuristic policies, and have considerably higher computational requirements than the PEIP policy.
Furthermore, these policies have more parameters and the BWB policy is arguably difficult to communicate to practitioners.

\begin{table*}[!htbp]
  \centering
  \caption{Summary of full factorial experiment}
	\fontsize{8pt}{9pt}\selectfont
  \label{tab:resultsOverall}
  \begin{tabular}{ l l l l l}
    \toprule
     Heuristic & Avg & Max & Min & \% \\
		\midrule
		VBS & 0.15 & 3.15 & -0.05 & 90 \\  
BWB & 0.03 & 2.67 & -0.29 & 42 \\  
DI & 0.71 & 4.97 & 0.00 & 100 \\   
CDI & 0.19 & 4.67 & -0.08 & 94 \\  
TBS & 55.11 & 373.88 & 1.68 & 100 \\
    \bottomrule
  \end{tabular}
\end{table*}

\section{Concluding Remarks}
\label{sec:concludingRemarks}
In this paper, we have proposed the PEIP policy for dual-sourcing inventory systems. This policy
uses an order-up-to rule for the expedited supplier, and places regular orders such that the expected expedited inventory position at the time of the regular order entering the expedited inventory position is raised to a target
level. In doing so, the PEIP policy leverages all the information contained in the pipeline -- as does the optimal policy.
This is in contrast with existing heuristic policies that either aggregate all outstanding orders or ignore them entirely and place a constant order every period regardless. We have further established a separability result that reduces the optimization of the PEIP policy to two simple one-dimensional optimization problems. 
An extensive numerical investigation has indicated that the PEIP policy outperforms the current best performing heuristic policies, even the ones that employ more control parameters than the PEIP policy.
The PEIP policy is thus vastly more advanced than existing heuristic policies, and yet retains an intuitive appeal for practitioners, has superior performance, and can be optimized efficiently.

We have shown that the PEIP policy is asymptotically optimal in both the long lead time and high service regimes. These results are useful for practice as it implies that the PEIP policy has asymptotic performance guarantees precisely when computation of optimal policies is prohibitive. 
Our proof techniques also extend beyond the dual sourcing system. The approach used for the long lead time regime also works for projection policies in lost sales inventory systems and we established high service asymptotic optimality of several other heuristic policies for the dual sourcing problem.

%Our numerical investigation has provided evidence that the PEIP policy converges to a TBS policy Furthermore, an interesting
%direction for future research would be to investigate whether the POL policy is also asymptotically
%optimal when the lead time difference between the regular supplier and the expedited supplier grows
%large. Exploiting the relation between the PEIP policy and the PIL policy for lost sales inventory systems
%would then be a promising approach. The latter policy has been shown to be asymptotically optimal for
%lost sales inventory systems when the lead time approaches infinity, albeit only under the assumption
%that demand is exponentially distributed.

% Acknowledgments here
\ACKNOWLEDGMENT{
An earlier version of this paper previously circulated under the title \textit{Effective dual-sourcing through inventory projection.} The research of the first author (PM) is supported by a starter grant from the Netherlands Organisation for Scientific Research (NWO), while the second author (MD) acknowledges support from NWO through grant no. VI.Veni.231E.060, and the research of the fourth author (CD) is supported by NWO through VENI grant no. VI.Veni.241E.058.

The authors acknowledge the use of ChatGPT Astra and Sol Models to verify and streamline select technical steps in our proofs. The overall arguments were constructed by the authors and they take full responsibility for all results in the paper. The authors also acknowledge useful discussions with Willem van Jaarsveld.
}

\bibliographystyle{informs2014}
\bibliography{bib}

\newpage
\begin{APPENDICES}

\section{Proof of Proposition \ref{prop:attainment}}
\label{sec:attainmentProof}

\begin{proof}{Proof of Proposition 1}
Fix $\Xb_t^r$ and an order $Q_t^r\geq0$ attaining $\E[O_{t+l}\mid X_t^r]=V$. Set $Q_{t+1}^r=0$,
and define the counterfactual process
\[
\widehat{O}_0=O_t+Q_{t-l+1}^r,\qquad
\widehat{O}_k=
(\widehat{O}_{k-1}+Q_{t-l+k+1}^r-D_{t+k-1})^+,
\quad k\in\{1,\ldots,l\}.
\]
Since $O_{t+1}\leq\widehat{O}_0$, monotonicity of the overshoot recursion and the i.i.d.\ demand assumption imply that $\E[\widehat{O}_l\mid \Xb_t^r]$ bounds the projected overshoot before regular ordering in period $t+1$, for every realization of $D_t$.

We show that $\widehat{O}_k\leq O_{t+k}+Q_{t-l+k+1}^r$.
Equality holds for $k=0$. By Equation (3) and $(x+q)^+\leq x^+ +q$ for $q\geq0$, induction gives
\[
\widehat{O}_k
\leq (O_{t+k-1}+Q_{t-l+k}^r
       +Q_{t-l+k+1}^r-D_{t+k-1})^+
\leq O_{t+k}+Q_{t-l+k+1}^r.
\]
Because we set $Q_{t+1}^r=0$, it follows that
\[
\E[\widehat{O}_l\mid \Xb_t^r]
\leq \E[O_{t+l}\mid \Xb_t^r]=V,
\]
that is, a zero order in the next period ensures we cannot exceed the projected overshoot level $V$.
The projected overshoot is continuous, non-decreasing,
and unbounded in the regular order size. Hence a
non-negative order $Q^r_{t+1}$ attains $V$ in period $t+1$.
Initial attainability and induction complete the proof. \hfill\Halmos
\end{proof}

\section{Proof of Lemma~\ref{lemma:stationarity}}
\label{section:StationarityProof}

In this appendix, we establish that, under the PEIP policy and a fixed lead time difference $l$, the Markov chain $\Xb^r$ converges to a unique steady state distribution.
The Markov chain $\Xb^r$ has transition probability kernel $
K(\xb,A):= \P \left\{\Xb^r_t \in A \mid \Xb^r_{t-1}=\xb \right\},
$
where $A\in\mathcal{B}\left(\R^l \right)$, the Borel $\sigma-$field on $\R^l$.
Then, $K^t(\cdot, \cdot)$ denotes the $t-$step transition probability kernel of $\Xb^r$.

The proof proceeds in three steps: (1) We establish the existence of an invariant probability measure for $\Xb^r$, $\pi$; (2) We use a coupling argument to prove that the invariant probability measure is unique; (3) We show that $K^t(\xb,\cdot)\to\pi(\cdot)$ with respect to a Wasserstein distance function for all attainable initial states $\Xb^r$.
These three steps form the subsections. Prior to this proof, we establish that the state space under the PEIP policy is compact.

Our analysis relies on the attainment result in Proposition~\ref{prop:attainment}, which ensures that the target level $V$ remains attainable in every period whenever it is attainable initially.
%We present a general state of $\Xb^r$ by $\xb=(o,q^r_{-1},\ldots,q^r_{-l+1})$.
%When $l=1$, $\xb=o$ and the pipeline is empty.
%When $l \geq 2$ and the system is in state $\xb$
%the overshoot $l$ periods later is defined as 
%\[
%    \OT_l(\xb,y):= \left(\left( \left(\left(\left( o+q^r_{-l+1}-D_{1}\right)^+ + q^r_{-l+2} - D_{2} \right)^+ + \dots \right)^+ +q^r_{-1}-D_{l-1} \right)^+ + y - D_l\right)^+.
%\]
% For $l=1$, only $\widetilde O_0(x)=o$ is needed.
%Observe that $\E \left[\OT_l(\xb,y) \right]$ is non-decreasing in $y$ and is unbounded on $\R_+$.
%Thus, by Proposition~\ref{prop:attainment} $\E \left[\OT_l(\xb,y) \right]=V$ has a unique solution denoted by $P_V(\xb)$.
Notice that
\[
    P_V(\xb) - \mu \leq \E \left[\OT_l(\xb,P_V(\xb)) \right]=V.
\]
Thus, $P_V(\xb) \leq V+\mu$, which implies that $Q^r_t\in[0,V+\mu]$.
Let $\XC$ be the state space of the Markov chain $\Xb^r$.
Under initial attainability, it can be expressed as
\begin{equation}
    \label{eq:StateSpace}
    \XC :=\left\{\xb \in \R_+^l: \E \left[\OT_l \left(\xb,0 \right) \right]  \leq V\right\}.
\end{equation}
Observe that $\XC$ is closed and bounded, and therefore compact.

\subsection{Existence of an invariant probability measure}

We will use the following definitions to prove that $\Xb^r$ has an invariant probability measure. 

\begin{definition}
    \label{def:weakFeller}
    The chain $\Xb^r$ is called \textit{weak Feller} if $K(\cdot, A)$ is a lower semicontinuous function for any open set $A \in \mathcal{B}(\R^l)$.
\end{definition}
\begin{definition}
    \label{def:boundedProb}
    The chain $\Xb^r$ is called \textit{bounded in probability on average} if for every $\epsilon>0$ there exists a compact subset $A \subset \XC$ such that $\displaystyle \liminf_{T \to \infty} \frac{1}{T} \sum_{t=1}^T K^t(\xb,A) \geq 1- \epsilon$, for each initial condition $\xb \in \XC$.
\end{definition}
Note that $\Xb^r$ is bounded in probability on average as
$\XC$ is compact.
\begin{lemma}
    \label{lemma:invariant_existance}
    \citep[Theorem 12.1.2][]{MT2009}
    If $\Xb^r$ is weak Feller and bounded in probability on average, then it admits at least one invariant probability measure.
\end{lemma}
Define the one-step shift operator, for $l\geq2$, by
\[
 \theta_d(\xb)=\left((o+q^r_{-l+1}-d)^+,q^r_{-l+2}, \dots,q^r_{-1}, 
            P_V(\xb)\right).
\]
% The coordinates following $P_V(\xb)$ are omitted when $l=2$.
When $l=1$, define $\theta_d(o)=(o+P_V(o)-d)^+$.
In what follows, we show that $\Xb^r$ is weak Feller. Together with Lemma~\ref{lemma:invariant_existance} and the compactness of the state space, this result ensures the existence of an invariant probability
measure $\pi$.

\begin{lemma}\label{lemma:weakFeller}
Under initial attainability, $\Xb^r$ is weak Feller.
\end{lemma}

\proof{Proof of Lemma~\ref{lemma:weakFeller}.}
The Markov chain $\Xb^r$ is weak Feller if and only if its transition kernel $K$ maps every bounded continuous function to a bounded continuous
function \citep[cf. Proposition 6.1.1][]{MT2009}. That is, for every bounded continuous function $f:\R_+^l\to\mathbb R$,
\[
    \int_{\R_+^l} f(\yb)K(\xb,d\yb)
    =\int_{\mathbb R_+}f(\theta_y(\xb))F_D(dy)
\]
is also bounded continuous on $\XC$, where $F_D$ is the demand distribution.
We restrict our argument to $\XC$ since the chain does not leave this set.

For any demand sample path $(d_1, \dots, d_l) \in \R_+^l$, note that $\OT_l(\xb,z)\mid D_1=d_1, \dots,D_l=d_l$ is continuous and uniformly bounded on $\XC \times [0,V+\mu]$.
Hence, by the bounded convergence theorem \citep[cf. Theorem 16.5][]{Billingsley1995} $\E \left[ \OT_l(\xb,z) \right]$ is jointly continuous in $\xb$ and $z$.
It follows that $P_V(\xb)$, $\theta_y(\xb)$, and $f \left(\theta_y(\xb) \right)$ are bounded and continuous in $\xb$ for every $y\geq0$.
Applying the bounded convergence theorem once more, we obtain that $\displaystyle \int_{\R_+}f \left(\theta_y(\xb) \right)F_D(dy)$ is bounded and uniformly continuous on $\XC$.
% For any demand sample path $(d_0,\ldots,d_{l-1})$, the random variable
% inside the expectation in \eqref{peipstat:forecast} is continuous in
% $(x,q)$ and uniformly bounded for $x\in\XC$ and
% $q\in[0,V+\mu]$. Hence, by the bounded convergence theorem
% (cf. Theorem 16.5 of Billingsley 1995), $F(x,q)$ is continuous in
% $(x,q)$. We must also verify continuity of its implicitly defined root.
% Let $x_n\to x$ in $\XC$. By
% \eqref{peipstat:order-bound}, every subsequence of $P_V(x_n)$ has a
% further convergent subsequence, with limit $q\in[0,V+\mu]$.
% Continuity gives $F(x,q)=V$, so uniqueness of the root implies
% $q=P_V(x)$. Thus $P_V(x_n)\to P_V(x)$.
% It follows that $P_V(x)$, $\theta_d(x)$, and $f(\theta_d(x))$ are
% bounded and continuous in $x$ for every $d\geq0$. Applying the bounded
% convergence theorem once more, we obtain that
% $\int_{\mathbb R_+}f(\theta_d(x))F_D(dd)$ is bounded and continuous on
% $\XC$.
\Halmos
\endproof

\begin{lemma}\label{lemma:pi_Existence}
The Markov chain $\Xb^r$ admits an invariant probability measure $\pi$.
\end{lemma}

\proof{Proof of Lemma~\ref{lemma:pi_Existence}.}
By Lemma~\ref{lemma:weakFeller}, $\Xb^r$ is weak Feller. It is also bounded in probability on average because $\XC$ is compact.
Therefore, by Lemma~\ref{lemma:invariant_existance}, an invariant probability measure $\pi$ exists.
\Halmos
\endproof

\subsection{Uniqueness of the invariant probability measure}

Next, we show that the invariant probability measure $\pi$ is unique.
The key idea of our argument is to couple two chains initialized from different stationary marginal distributions and driven by the same demand sequence.
That is, let $\XbT^r=\{(\Xb^r_{1,t},\Xb^r_{2,t})\}_{t\geq0}$ be a tuple of two PEIP systems driven by the exact same demand process, but the initial distributions of $\Xb^r_{1,0}$ and $\Xb^r_{2,0}$ can be distinct.
Both systems use the same projected overshoot level $V=\xi-S^e$. Let $ \Tilde{K}$ denote the transition probability kernel of $\XbT^r$ and $\Tilde{\XC}:=\XC\times\XC$ its state space.
Our first result is that $\XbT^r$ has an invariant probability measure in which the marginal measures of $\Xb^r_1$ and $\Xb^r_2$ are also invariant probability measures of $\Xb^r$.

\begin{lemma}\label{lemma:XtildeInvariantMeasure}
For any two invariant probability measures $\pi_1$ and $\pi_2$
of $\Xb^r$, the coupled process $\XbT^r$ has an invariant
probability measure with marginals $\pi_1$ and $\pi_2$.
\end{lemma}

\proof{Proof of Lemma~\ref{lemma:XtildeInvariantMeasure}.}
Lemma~\ref{lemma:weakFeller} applies to $\Xb^r_i$, $i\in\{1,2\}$ and so $\XbT^r$ is weak Feller.
Its state space $\Tilde{\XC}$ is compact because $\XC$ is compact, and so by Lemma~\ref{lemma:pi_Existence}, $\XbT^r$ has an invariant probability measure.
 We now construct such a invariant probability measure where the marginal measures of $\Xb^r_1$ and $\Xb^r_2$ are invariant for $\Xb^r$ by using the Krylov-Bogoliubov Construction.

Let $\pi_1$ and $\pi_2$ be two (possibly distinct) invariant probability measures of $\Xb^r$, which exist by Lemma~\ref{lemma:pi_Existence}. 
Let the coupled process be initialized with a joint distribution $\nu_0$ whose marginals are $\pi_1$ and $\pi_2$ (e.g., the product measure $\pi_1\times\pi_2$).
Let $\nu_t$ denote the joint distribution after $t$ transitions.
Because $\pi_1$ and $\pi_2$ are individually invariant, the marginal distributions of $\nu_t$ remain $\pi_1$ and $\pi_2$ for all $t\geq0$.
For each $t\geq1$, define the time-averaged probability measure
\begin{equation}
    \label{eq:Def_averagedProb}
    \nubar_t(A):=\frac{1}{t}\sum_{n=1}^t\nu_n(A)
    =\frac{1}{t}\sum_{n=1}^t \int_{\Tilde{\XC}} \Tilde{K}^n(\xb,A)\nu_0(d\xb).
\end{equation}
Because $\Tilde{\XC}$ is compact, there exists a subsequence of $\{\nubar_t\}$, indexed by some $\mathcal T\subseteq\N$, and a probability measure $\Tilde{\pi}$ such that, as $t\to\infty$ with $t\in\mathcal T$, $\nubar_t$ converges weakly to $\Tilde{\pi}$ by Prohorov's theorem \citep[cf. Theorem 5.1][]{Billingsley1999}.

For any bounded continuous function $f$ on $\Tilde{\XC}$, define
\[
  \Tilde{K}f(\xbT)
   =\int_{\Tilde{\XC}}
       f(\ybT) \Tilde{K}(\xbT,d\ybT).
\]
The telescoping sum in \eqref{eq:Def_averagedProb} ensures that
\[
    \Bigg\vert \int_{\Tilde{\XC}}\Tilde{K}f(\xbT)\nubar_t (d\xbT) - \int_{\Tilde{\XC}} f(\xbT) \nubar_t (d\xbT) \Bigg\vert \leq \frac{2 \sup_{\xbT \in \Tilde{\XC}} \vert f(\xbT) \vert}{t},
\]
which vanishes as $t\to\infty$. Because the coupled chain is weak
Feller, $ \Tilde{K}f$ is also a bounded continuous function.
Weak convergence therefore implies that
\[
    \int_{\Tilde{\XC}} \Tilde{K}f(\xbT)\nubar_t (d\xbT) \to \int_{\Tilde{\XC}} \Tilde{K}f(\xbT) \Tilde{\pi} (d\xbT),
     \quad \mbox{and,} \quad 
     \int_{\Tilde{\XC}} f(\xbT) \nubar_t (d\xbT) \to \int_{\Tilde{\XC}} f(\xbT) \Tilde{\pi} (d\xbT).
\]
Comparing these limits yields 
\[
    \int_{\Tilde{\XC}} \Tilde{K}f(\xbT) \Tilde{\pi} (d\xbT) = \int_{\Tilde{\XC}} f(\xbT) \Tilde{\pi} (d\xbT).
\]
Hence, $\Tilde{\pi}$ is an invariant probability measure for the
coupled process. Furthermore, because every $\nu_n$ has marginals
$\pi_1$ and $\pi_2$, the limit measure inherits exactly the
same marginals.\hfill\Halmos
\endproof

Next we show that the coupled process $\XbT^r$ satisfies a
useful contraction inequality. For any variable $Y$ in the coupled
system, let $Y_{1,t}$ and $Y_{2,t}$ denote its value in the first and
second system, respectively, and define the difference as
$\Delta Y_t:=Y_{1,t}-Y_{2,t}$. Since the state space
$\Tilde{\XC}$ is compact, $Q^r_{i,t}$ and $O_{i,t}$
are bounded and possess finite moments of all orders.

\begin{lemma}
\label{prop:contracioniIneqCoupledProcess}
The coupled process satisfies the following contraction inequalities for all $t\geq l-1$:
\begin{equation}
    \label{eq:contraction-eq}
    \E \left[(\Delta O_{t+1})^2 \right] \leq\E \left[(\Delta O_t)^2 \right]-\E \left[(\Delta O_{t+1}-\Delta O_t)^2 \right],
\end{equation}
\begin{equation}
    \label{eq:order-identification}
    V\E \left[\big \vert\Delta Q^r_{t-l+1} \big \vert \right] \leq (V+l \mu)\left(\E \left[\big \vert\Delta O_t \big \vert \right] +\E \left[\big \vert \Delta O_{t+1} \big \vert \right]\right).
\end{equation}
\end{lemma}

\proof{Proof of Lemma~\ref{prop:contracioniIneqCoupledProcess}.}
We establish a set of path-wise identities for $t\geq l-1$.
Expanding $\E \left[(\Delta O_{t+1}-\Delta O_t)^2 \right]$ yields
\begin{equation}
    \label{eq:quadratic_expansion}
    \E \left[(\Delta O_{t+1}-\Delta O_t)^2 \right] = \E \left[(\Delta O_{t+1})^2 \right] + \E \left[(\Delta O_t)^2\right] - 2 \E \left[\Delta O_{t+1}\Delta O_t \right].
\end{equation}
Recall the overshoot recursion in Equation~\eqref{eq:Dynamics_O}:
\[
 O_{i,t+1}= \left(O_{i,t}+Q^r_{i,t-l+1}-D_t \right)^+,\qquad i\in\{1,2\}.
\]
Observe that for any real numbers $a$ and $b$
%that $(\cdot)^+$ is (1-)lipschitz and monotone, therefore, for any real numbers $a$ and $b$,
\[
    (a^+-b^+)^2 =\vert a^+-b^+  \vert^2 \leq \vert a^+-b^+  \vert \vert a-b  \vert = \left(a^+-b^+\right)(a-b).
\]
This relation, together with  Equation~\eqref{eq:Dynamics_O}, gives
\begin{equation}
    \label{eq:lipschitzDeltaQ}
    \left(\Delta O_{t+1} \right)^2 \leq \Delta O_{t+1} (\Delta O_t+\Delta Q_{t-l+1}^r).
\end{equation}
By the independence of future demand from the current state,
\begin{equation*}
    \label{eq:zero-conditional-mean}
    \E \left[O_{i,t+1} \big \vert \XbT_{t-l+1}^r \right]=V\qquad i\in\{1,2\} \qquad
    \Rightarrow
    \E \left[\Delta O_{t+1} \big \vert \XbT_{t-l+1}^r \right]=0.
\end{equation*}
The order decision difference $\Delta Q_t^r$ is known conditional on
$\XbT_t^r$. When substituting \eqref{eq:lipschitzDeltaQ} into \eqref{eq:quadratic_expansion} we need to deal with a cross term by conditioning on $\XbT_{t-l+1}^r$ as follows
\begin{equation}
    \label{eq:cross-term}
    \E \left[\Delta O_{t+1}\Delta Q_{t-l+1}^r \right] = \E\left[\Delta Q_{t-l+1}^r \E \left[\Delta O_{t+1}\big \vert \XbT_{t-l+1}^r \right]\right]=0.
\end{equation}
Combining \eqref{eq:quadratic_expansion}, \eqref{eq:lipschitzDeltaQ}, and \eqref{eq:cross-term} gives \eqref{eq:contraction-eq}.

It remains to relate regular-order differences to overshoot
differences. Write $b=\Delta Q_{t-l+1}^r$. We claim the pathwise inequality
\begin{equation}
    \label{eq:order-pathwise}
    \big \vert \Delta Q_{t-l+1}^r \big \vert\left(O_{1,t+1}\ind_{\{\Delta Q_{t-l+1}^r>0\}} +O_{2,t+1}\ind_{\{\Delta Q_{t-l+1}^r \leq 0\}}\right) \leq (V+l \mu)\left(\big \vert\Delta O_t\big \vert+\big \vert\Delta O_{t+1}\big \vert \right).
\end{equation}
If the overshoot selected on the left hand side of the inequality is zero, the inequality is immediate.
If it is positive and the other overshoot is zero, the selected overshoot is equal to $\big \vert \Delta O_{t+1}\big \vert$, and $\big \vert \Delta Q_{t-l+1}^r \big \vert \leq V+l \mu$.
If both overshoots are positive, then by \eqref{eq:Dynamics_O}
$\Delta O_{t+1}=\Delta O_t+ \Delta Q_{t-l+1}^r$ and each overshoot is at most $V+l \mu$, which again gives \eqref{eq:order-pathwise}.
Taking expectations in \eqref{eq:order-pathwise} and conditioning
its left hand side on $\XbT_{t-l+1}^r$,
\[
    \E \left[\E \left[ \big \vert \Delta Q_{t-l+1}^r \big \vert\left(O_{1,t+1} \ind_{\{\Delta Q_{t-l+1}^r>0\}} +O_{2,t+1} \ind_{\{\Delta Q_{t-l+1}^r\leq 0\}}\right) \Big \vert \XbT_{t-l+1}^r\right] \right] = V \E \left[ \big \vert \Delta Q_{t-l+1}^r \big \vert \right],
\]
which yields \eqref{eq:order-identification}.
\hfill\Halmos
\endproof

Now we can use this contraction result together with the existence of
an invariant probability measure of $\XbT^r$ to show that
there is a unique invariant probability measure of $\Xb^r$.

\begin{proposition}\label{proposition:uniqueness}
The Markov chain $\Xb^r$ has a unique invariant probability measure.
\end{proposition}

\proof{Proof of Proposition~\ref{proposition:uniqueness}.}
Assume for the sake of contradiction that the invariant probability measure of $\Xb^r$ is not unique and let $\pi_1$ and $\pi_2$ be two distinct invariant probability measures of the Markov chain $\Xb^r$.
Let $\Tilde{\pi}$ be an invariant probability measure for $\tilde\Xb^r$ such that its marginals are given by $\pi_1$ and $\pi_2$, respectively.
These measures exist by Lemmas~\ref{lemma:invariant_existance} and~\ref{lemma:XtildeInvariantMeasure}.

Because the coupled chain operates strictly in stationarity under $\Tilde{\pi}$, the expected squared differences must be constant over time, that is, $\E[(\Delta O_{t+1})^2]=\E[(\Delta O_t)^2]$.
Combining this stationarity condition with \eqref{eq:contraction-eq} forces $\Delta O_{t+1}=\Delta O_t$ almost surely. Stationarity extends this equality to every $t\geq0$. Consequently, $\Delta O_{t+1}=\Delta O_{t-l+1}$ almost surely for $t\geq l-1$.
Hence, $\Delta O_{t-l+1} = \E[\Delta O_{t-l+1}\mid\XbT^r_{t-l+1}] = \E[\Delta O_{t+1}\mid\XbT^r_{t-l+1}] = 0$, so $\Delta O_t=0$ a.s.
% The latter variable is measurable with respect to $\XbT^r_{t-l+1}$, whereas \eqref{peipstat:zero-conditional-mean} gives
% \[
%  0=\E[\Delta O_{t+1}\mid\widetilde X^r_{t-l+1}]
%    =\Delta O_{t-l+1}\quad\hbox{almost surely}.
% \]
Inequality \eqref{eq:order-identification}, together with $V>0$, now forces $\Delta Q_t^r=0$ almost surely for every $t\geq0$.
After shifting through the pipeline, all state coordinates coincide; by stationarity, this is also true of the initial coupled state.
Hence $\Xb_1^r=\Xb_2^r$, $\Tilde{\pi}$-almost surely, which implies that their marginal distributions are identical, proving $\pi_1=\pi_2$. But this implies that invariant probability measures of $\Xb^r$ cannot be distinct and the invariant probability measure is unique.\hfill\Halmos
\endproof

The following corollary follows directly from the proof of Proposition~\ref{proposition:uniqueness} and we will use the result later.

\begin{corollary}
    \label{corollary:X1X2as}
    Every invariant probability measure of $\XbT^r$, $\Tilde{\pi}$, is such that $\Xb^r_1 = \Xb^r_2$, $\Tilde{\pi}$-almost surely.
\end{corollary}

\subsection{Convergence to the invariant probability measure}

Finally, we show that the probability measure of $\Xb_t^r$, i.e., $K^t(x,\cdot)$, converges weakly to the unique invariant probability measure $\pi$ as $t\to\infty$, regardless of the attainable initial state $\Xb^r_0$. We couple two chains operating under the PEIP policy: one initialized from an arbitrary state in $\XC$, and the other initialized in stationarity.
\begin{definition}
    \label{def:Wasserstein}
    The \textit{$a-$Wasserstein} distance between two probability measures $\eta_1$ and $\eta_2$ is defined by
\begin{equation}
    \label{eq:Wasserstein}
    W_a \left(\eta_1,\eta_2 \right):= \inf_{Y_1 \sim \eta_1,Y_2 \sim\eta_2} \left(\E \left[ \Vert Y_1 - Y_2\Vert_a^a\right] \right)^{\frac{1}{a}}, \qquad a\geq1.
\end{equation}
\end{definition}
We use the $a$-Wasserstein distance to show that the distance between the marginal distributions converges to zero as $t\to\infty$. All probability measures on the compact state space $\XC$ have finite moments of every order, so this distance is well defined for every finite $a\geq1$.

\begin{proposition}\label{prop:pi_convergence}
The probability measure of $\Xb_t^r$ converges to $\pi$, for any initial state $\xb\in\XC$. That is, for every $\xb\in\XC$ and finite $a\geq1$,
\[
 \lim_{t\to\infty}W_a(K^t(\xb,\cdot),\pi)=0.
\]
\end{proposition}

\proof{Proof of Proposition~\ref{prop:pi_convergence}.}
Consider again $\tilde\Xb^r$ and suppose $\Xb^r_{1,0}=\xb$, where $\xb$ is an arbitrary attainable initial state, and $\Xb^r_{2,0}$ is distributed as the invariant probability measure of $\Xb^r$, $\pi$.
By \eqref{eq:contraction-eq}, $\E[(\Delta O_t)^2]$ is non-increasing for $t\geq l-1$ and
converges to a constant $a_\infty\geq0$. It remains to show that
$a_\infty=0$.

Let $\nubar_t$ denote the averaged probability measure after $t$ transitions, as defined in \eqref{eq:Def_averagedProb}.
Then, for every bounded continuous function $f$ on $\Tilde{\XC}$,
\begin{equation}
    \label{eq:averagedFunc}
    \int_{\Tilde{\XC}} f(\yb) \nubar_t(d\yb) = \frac{1}{t} \sum_{n=1}^t \E \left[ f \left(\XbT^r_n \right) \right].
\end{equation}
By compactness of $\Tilde{\XC}$, the weak Feller property, and the weak convergence of $\{\nubar_t\}$ to $\Tilde{\pi}$, it follows that as $t \to \infty$ with $t \in \mathcal{T} \subseteq \N_0$
\begin{equation}
    \label{eq:convergenceFunc}
    \int_{\Tilde{\XC}} f(\yb) \nubar_t(d\yb) - \int_{\Tilde{\XC}} f(\yb) \Tilde{\pi}(d\yb) \to 0.
\end{equation}

% Compactness supplies a weakly convergent subsequence
% $\nubar_t\Rightarrow\Tilde{\pi}$, $t\in\mathcal T$.
% The weak Feller property and the telescoping bound
% \eqref{peipstat:telescoping} show that its limit is invariant, exactly
% as in the proof of Lemma~\ref{peipstat:stationary-coupling}; that part
% of the proof does not require stationary initial marginals.

\noindent
Choose $f(y)=(o_1-o_2)^2$, which is continuous and bounded on $\Tilde{\XC}$.
By \eqref{eq:Def_averagedProb}
\[
    \int_{\Tilde{\XC}} f(\yb) \nubar_t(d\yb) - \frac{1}{t} \sum_{n=1}^t \E \left[\left(\Delta O_n \right)^2 \right] = 0 \quad \Rightarrow \quad \int_{\Tilde{\XC}} f(\yb) \nubar_t(d\yb) \to a_\infty.
\]
On the other hand, along the weakly convergent subsequence,
\[
    \int_{\Tilde{\XC}} f(\yb) \nubar_t(d\yb) - \int_{\Tilde{\XC}} f(\yb) \Tilde{\pi}(d\yb) \to 0 \quad \Rightarrow \quad \int_{\Tilde{\XC}} f(\yb) \nubar_t(d\yb) \to 0.
\]
The last equality follows because, under any invariant coupling, $\Delta O_t=0$ almost surely by Corollary~\ref{corollary:X1X2as}.
Therefore $a_\infty=0$, i.e., $\lim_{t\to\infty}\E[(\Delta O_t)^2]=0$.
Together with \eqref{eq:order-identification}, this result implies $\E[\big\vert \Delta Q_t^r \big \vert ]\to0$.
Since $\big \vert \Delta Q_t^r\big \vert \leq V+\mu$, it also implies $\E[(\Delta Q_t^r)^2]\to 0$.
Moreover,
\begin{equation}
    \label{eq:2norm}
    \E\left[\big\Vert \Xb^r_{1,t}-\Xb^r_{2,t}\big\Vert _2^2 \right]
    =\E \left[(\Delta O_t)^2 \right] +\sum_{j=1}^{l-1}\E \left[(\Delta Q^r_{t-j})^2 \right].
\end{equation}
where $\Vert \cdot\Vert_2$ is the $2-$norm.
All the terms on the right hand side of \eqref{eq:2norm} converge to zero, hence 
\[
    \lim_{t\to \infty} \E \left[\Vert \Xb^r_{1,t} - \Xb^r_{2,t} \Vert_2^2 \right]=0.
\]
Because $\XC$ is a compact set in Euclidean space, for any finite $a\geq1$,
\[
    \lim_{t\to \infty} \E \left[\Vert \Xb^r_{1,t} - \Xb^r_{2,t}\Vert_a^a \right] =0.
\]
Observe that $\Xb^r_{1,t}$ is distributed according to $K^t(\xb^r,\cdot)$
and $\Xb^r_{2,t}$ according to $\pi$. It follows from the definition of $W_a$ that
\[
    \lim_{t \to \infty}W_a \left( K^t(\xb,\cdot) , \pi\right) \leq  \lim_{t\to \infty} \E \left[\Vert \Xb^r_{1,t} - \Xb^r_{2,t}\Vert_a^a \right] ^{1/a}=0.
\]
Together with Proposition~\ref{proposition:uniqueness}, this proves Lemma~\ref{lemma:stationarity}.\hfill\Halmos
\endproof

\section{Detailed results numerical investigation}
\label{sec:detailedResults}
Table \ref{tab:resultsTestbedFull} presents the percentage gap in the cost-rate of each benchmark heuristic policy relative to the PEIP policy. In this table, we first distinguish between subsets of instances with the
same value for a specific input parameter of Table \ref{tab:testbed1} and then present the results for all instances. 

\begin{table*}[!htb]
  \centering
  \caption{Detailed result full factorial experiment}
	\fontsize{8pt}{9pt}\selectfont
	\setlength{\tabcolsep}{2pt}
  \label{tab:resultsTestbedFull}
  \begin{tabularx}{1\textwidth}{X c l l l l c@{\hspace{0.0cm}} l l l l c@{\hspace{0.0cm}} l l l l c@{\hspace{0.0cm}} l l l l c@{\hspace{0.0cm}} l l l l }
  \toprule
	&  & \multicolumn{4}{c}{VBS} & & \multicolumn{4}{c}{BWB} & & \multicolumn{4}{c}{DI} & & \multicolumn{4}{c}{CDI} & & \multicolumn{4}{c}{TBS} \\
	\cmidrule{3-6}\cmidrule{8-11}\cmidrule{13-16} \cmidrule{18-21} \cmidrule{23-26}
	Parameter &Value & Avg & Max & Min & \% & & Avg & Max & Min & \% & & Avg & Max & Min & \% & & Avg & Max & Min & \%   & & Avg & Max & Min & \% \\
		\midrule
		\multirow{ 5}{2.5cm}{$p$}
&	4	& 0.32 & 0.65 & 0.08 & 100 && 0.09 & 0.45 & -0.11 & 78 && 1.76 & 4.97 & 0.23 & 100 && 0.36 & 4.19 & 0.06 & 100 && 10.73 & 30.99 & 1.68 & 100 \\
&	9	& 0.20 & 0.62 & 0.02 & 100 && 0.01 & 0.27 & -0.21 & 47 && 1.01 & 3.27 & 0.11 & 100 && 0.28 & 4.32 & 0.02 & 100 && 21.74 & 58.51 & 4.48 & 100 \\
&	19	& 0.08 & 0.33 & -0.05 & 96 && -0.04 & 0.66 & -0.29 & 22 && 0.50 & 1.73 & 0.04 & 100 && 0.15 & 3.10 & 0.01 & 100 && 38.35 & 90.34 & 9.34 & 100 \\
&	49	& 0.12 & 3.15 & -0.03 & 83 && -0.03 & 0.46 & -0.18 & 26 && 0.17 & 0.61 & 0.00 & 100 && 0.04 & 0.43 & -0.07 & 92 && 76.82 & 199.62 & 23.71 & 100 \\
&	99	& 0.02 & 0.09 & -0.02 & 72 && 0.12 & 2.67 & -0.13 & 35 && 0.11 & 2.81 & 0.00 & 100 && 0.12 & 4.67 & -0.08 & 79 && 127.90 & 373.88 & 32.34 & 100 \\
\addlinespace	\multirow{ 4}{2.5cm}{$l^r$}
&	2	&	0.07 & 0.29 & 0.00 & 94 && 0.08 & 2.67 & -0.05 & 51 && 0.27 & 1.36 & 0.00 & 100 && 0.05 & 0.29 & -0.01 & 90 && 48.75 & 235.21 & 1.75 & 100 \\
&	3	&	0.15 & 2.64 & -0.01 & 92 && 0.02 & 0.73 & -0.10 & 42 && 0.58 & 2.80 & 0.01 & 100 && 0.17 & 4.32 & -0.02 & 96 && 54.35 & 286.08 & 1.68 & 100 \\
&	4	&	0.15 & 0.54 & -0.01 & 90 && 0.02 & 0.85 & -0.18 & 38 && 0.85 & 3.86 & 0.00 & 100 && 0.26 & 4.67 & -0.07 & 97 && 57.45 & 334.64 & 1.76 & 100 \\
&	5	&	0.21 & 3.15 & -0.05 & 84 && 0.00 & 1.31 & -0.29 & 36 && 1.14 & 4.97 & 0.01 & 100 && 0.27 & 4.19 & -0.08 & 94 && 59.87 & 373.88 & 1.70 & 100 \\
\addlinespace	\multirow{ 6}{2.5cm}{$\frac{\sigma}{\mu}$}
&	0.15	&	0.24 & 3.15 & -0.02 & 97 && 0.09 & 2.67 & -0.14 & 52 && 0.81 & 4.97 & 0.00 & 100 && 0.28 & 4.67 & -0.08 & 92 && 76.50 & 373.88 & 4.29 & 100 \\
&	0.25	&	0.15 & 0.61 & 0.00 & 100 && 0.09 & 2.47 & -0.18 & 42 && 0.81 & 4.87 & 0.01 & 100 && 0.13 & 1.46 & -0.02 & 92 && 63.18 & 240.11 & 3.85 & 100 \\
&	0.5	    &	0.12 & 0.62 &  0.00 & 98 && -0.03 & 0.85 & -0.20 & 30 && 0.79 & 4.39 & 0.01 & 100 && 0.24 & 4.32 & -0.01 & 95 && 60.73 & 292.31 & 3.38 & 100 \\
&	1	    &	0.10 & 0.60 & -0.01 & 83 && -0.03 & 0.44 & -0.18 & 25 && 0.72 & 3.37 & 0.02 & 100 && 0.13 & 0.79 & -0.01 & 92 && 49.06 & 228.77 & 2.68 & 100 \\
&	1.5	    &	0.11 & 0.44 & -0.02 & 78 && 0.00 & 0.20 & -0.21 & 43 && 0.61 & 2.62 & 0.01 & 100 && 0.15 & 0.52 & 0.00 & 98 && 43.74 & 214.23 & 2.11 & 100 \\
&	2	    &   0.15 & 0.55 & -0.05 & 85 && 0.05 & 0.45 & -0.29 & 58 && 0.54 & 2.81 & 0.01 & 100 && 0.19 & 2.79 & 0.00 & 97 && 37.43 & 177.20 & 1.68 & 100 \\
\addlinespace	\multirow{ 3}{2.5cm}{$c^e=\delta p l$}
&	0.1	&	0.20 & 0.65 & -0.02 & 93 && 0.02 & 0.45 & -0.18 & 51 && 1.15 & 4.97 & 0.03 & 100 && 0.27 & 4.67 & -0.08 & 90 && 24.68 & 88.41 & 1.70 & 100 \\
&	0.2	&	0.15 & 0.62 & -0.03 & 89 && 0.02 & 2.67 & -0.20 & 37 && 0.78 & 3.39 & 0.02 & 100 && 0.21 & 4.19 & -0.01 & 94 && 48.60 & 182.03 & 4.48 & 100 \\
&	0.4	&	0.11 & 3.15 & -0.05 & 89 && 0.06 & 2.47 & -0.29 & 44 && 0.43 & 2.81 & 0.00 & 100 && 0.15 & 2.79 & -0.01 & 97 && 89.73 & 373.88 & 9.82 & 100 \\
\addlinespace 
Total & &  0.15 & 3.15 & -0.05 & 90 && 0.03 & 2.67 & -0.29 & 42 && 0.71 & 4.97 & 0.00 & 100 && 0.19 & 4.67 & -0.08 & 94 && 55.11 & 373.88 & 1.68 & 100 \\
  \bottomrule
  \end{tabularx}
\end{table*}

\end{APPENDICES}

\end{document}